\documentclass[11pt]{article}
\def\mytitle{An iterative method
using boundary distance
for box-constrained nonlinear semidefinite programs}

\usepackage{amssymb}
\usepackage{amsthm}
\usepackage{braket}
\usepackage{ascmac}

\usepackage{color}
\usepackage{enumerate}

\usepackage{atbegshi}
\usepackage[
dvipdfmx,
bookmarks=true,
bookmarksnumbered=true, 
bookmarkstype=toc,
setpagesize=false,
pdftitle={
\mytitle
},
pdfauthor={Akihiko Komatsu, Makoto Yamashita},
pdfkeywords={Nonlinear semidefinite programs, Global convergence
}
]{hyperref}

\def\@themcountersep{}
\newtheorem{THEO}{Theorem}[section]
\newtheorem{ALGO}[THEO]{Algorithm}

\newtheorem{CORO}[THEO]{Corollary}

\newtheorem{LEMM}[THEO]{Lemma}

\def\bold#1{\mbox{\boldmath $#1$}}

\def\0{\mbox{\bf 0}}
\def\1{\mbox{\bf 1}}
\def\2{\mbox{\bf 2}}
\def\3{\mbox{\bf 3}}
\def\4{\mbox{\bf 4}}
\def\5{\mbox{\bf 5}}
\def\6{\mbox{\bf 6}}
\def\7{\mbox{\bf 7}}
\def\8{\mbox{\bf 8}}
\def\9{\mbox{\bf 9}}

\newdimen\zhige \zhige=0pt
\def\chige#1{{\setbox\zhige\hbox{#1}\ifdim\ht\zhige=1ex\accent24 #1%
  \else\ooalign{\unhbox\zhige\crcr\hidewidth\char24\hidewidth}\fi}}

\def\l{\mbox{\boldmath $l$}}

\def\p{\mbox{\boldmath $p$}}

\def\u{\mbox{\boldmath $u$}}
\def\v{\mbox{\boldmath $v$}}

\def\x{\mbox{\boldmath $x$}}

\def\A{\mbox{\boldmath $A$}}
\def\B{\mbox{\boldmath $B$}}
\def\C{\mbox{\boldmath $C$}}
\def\D{\mbox{\boldmath $D$}}

\def\I{\mbox{\boldmath $I$}}

\def\K{\mbox{\boldmath $K$}}
\def\L{\mbox{\boldmath $L$}}
\def\M{\mbox{\boldmath $M$}}

\def\O{\mbox{\boldmath $O$}}
\def\P{\mbox{\boldmath $P$}}
\def\Q{\mbox{\boldmath $Q$}}

\def\S{\mbox{\boldmath $S$}}

\def\U{\mbox{\boldmath $U$}}
\def\V{\mbox{\boldmath $V$}}
\def\W{\mbox{\boldmath $W$}}
\def\X{\mbox{\boldmath $X$}}
\def\Y{\mbox{\boldmath $Y$}}

\def\FC{\mbox{$\cal F$}}

\def\KC{\mbox{$\cal K$}}
\def\LC{\mbox{$\cal L$}}

\def\OC{\mbox{$\cal O$}}

\def\Real{\mbox{$\mathbb{R}$}}

\def\SMAT{\mbox{$\mathbb{S}$}}

\ifx11 

\fi

\begin{document}

\noindent \textcolor{blue}{\leaders\hrule width 0pt height 0.08cm \hfill}
\vspace{0.5cm} 

{\Large \bf \noindent 
\mytitle
}

\vspace{0.2cm}
 \noindent
Akihiko Komatsu
\footnote{
Equities Department, Tokyo Stock Exchange, INC., 
2-1, Nihombashi-kabuto-cho, Chuo-ku, Tokyo 103-8220, Japan. 
}
and
Makoto Yamashita
\footnote{
(Corresponding Author) Department of Mathematical and Computing Sciences,
 Tokyo Institute of Technology, 2-12-1-W8-29 Ookayama, Meguro-ku, Tokyo
 152-8552, Japan (Makoto.Yamashita@is.titech.ac.jp).
 The work of the second author was supported by JSPS KAKENHI
Grant Numbers 24710151.
 }
\\
Submitted: May 14, 2015.

\vspace{0.3cm}

\noindent {\bf Abstract:} 

We propose an iterative method for nonlinear semidefinite programs
with box constraints.
The search direction in the proposed method utilizes the distance 
from the current point to the boundary of a feasible set.
The computation of the search direction exploits the second derivative of the objective function
only in a quadratic form, and this property saves the computation cost
compared to an evaluation of the whole entries of the second derivative.
We compute a step length in an interval determined by a radius
and we update
the radius using a quadratic approximation function.
In this paper,
we also discuss convergence properties of the proposed method 
based on structures of the search direction.
Numerical tests show that the proposed method solves problems
in which the size of a variable matrix is larger than 5,000
and that it is faster than a feasible direction method
for objective functions with strong nonlinearity.


\vspace{0.2cm}
\noindent {\bf Keywords:} \newline
Nonlinear programming, Semidefinite programming, Box constraint

\noindent {\bf AMS Classification:} \newline
90C22 Semidefinite programming, 
90C30 Nonlinear programming.

\noindent \textcolor{blue}{\leaders\hrule width 0pt height 0.04cm \hfill}

\section{Introduction}

This paper is concerned with  a box-constrained nonlinear semidefinite problem
(shortly, box-constrained SDP)
\begin{eqnarray}
\min \quad f(\X) \quad \mbox{subject to} \quad \O \preceq \X \preceq \I.
\label{eq:boxSDP}
\end{eqnarray}
The variable in this problem is $\X \in \SMAT^n$, and we use 
$\SMAT^n$ to denote the space of $n \times n$ symmetric matrices.
The notation $\A \succeq \B$ for $\A, \B \in \SMAT^n $ means that 
the matrix $\A-\B$ is positive semidefinite.
The matrix $\I$ is the identity matrix of the appropriate dimension.
We assume 
that the objective function $f : \SMAT^n \to \Real$ 
is a twice continuously differentiable function on an open set containing 
the feasible set $\FC := \{\X \in \SMAT^n : \O \preceq \X \preceq \I\}$.

The feasible set of (\ref{eq:boxSDP}) can express a more general feasible set
 $\{\X \in \SMAT^n : \L \preceq \X \preceq \U\}$ with $\L,\U \in \SMAT^n$ such that 
$\L \preceq \U$. 
This type of problems appears as a sub problem in other methods~\cite{JI14}.
We can assume that $\U-\L$ is positive definite
without loss of generality~\cite{XU11}, therefore,
we use a Cholesky factorization matrix $\C$ of $\U-\L$ that satisfies
$\U-\L = \C\C^T$ to convert a problem
\begin{eqnarray*}
\min \quad f(\X) \quad \mbox{subject to} \quad \L \preceq \X \preceq \U
\end{eqnarray*}
into an equivalent problem
\begin{eqnarray*}
\min \quad f(\C\overline{\X}\C^T+\L) \quad \mbox{subject to} \quad \O \preceq \overline{\X} \preceq \I
\end{eqnarray*}
by the relation $\overline{\X} = \C^{-1} (\X - \L) (\C^{-1})^T$. 
In this paper, we use the superscript $T$ to denote the transpose of a matrix.

A box-constrained nonlinear optimization problem
\begin{eqnarray}
\min \quad f(\x) \quad \mbox{subject to} \quad \l \le \x \le \u, \x \in \Real^n.
\label{eq:boxLP}
\end{eqnarray}
is an important case of (\ref{eq:boxSDP}), since
if the variable matrix $\X$ in (\ref{eq:boxSDP}) is a diagonal matrix, 
(\ref{eq:boxSDP}) is reduced to  (\ref{eq:boxLP}).
The problem (\ref{eq:boxLP}) is a basic problem in constrained optimization
and many methods are proposed. 
Hei {\it et. al.}~\cite{HEI08} compared  the performance of four active-set methods 
and two interior-point methods. 
Trust-region methods for (\ref{eq:boxLP}) are also discussed in \cite{COLEMAN96, CONN88, WANG13},

On the other hand, the positive semidefinite condition on a matrix ($\X \succeq \O$)
 is extensively studied in the context of SDP (semidefinite programs).
The range of SDP applications is very wide and includes control theory~\cite{BOYD94}, 
combinatorial optimization~\cite{GOEMANS95}, polynomial 
optimization~\cite{LASSERRE01} and quantum chemistry~\cite{FUKUDA07}. 
Many software packages, for example \cite{TODD99, YAMASHITA12}, have been developed for SDP.
A number of studies on SDP can be found in the survey of Todd~\cite{TODD01},
the handbook edited by Anjos and Lassere~\cite{ANJOS12} and the references therein.


For solving the box-constrained SDP~(\ref{eq:boxSDP}),
we may apply the penalty barrier method proposed in \cite{BENTAL97, KOCVARA03}.
Though it can handle the problem (\ref{eq:boxSDP}) with additional constraints,
it requires the full information of the second derivative of the objective function,
and it can solve the problems in practical time only when the size of variable matrix is small; $n \le 500$.

To solve large problems with $n \ge 500$, 
we should discuss methods specialized for solving 
(\ref{eq:boxSDP}).
Xu {\it et al}~\cite{XU11} proposed a feasible direction method for (\ref{eq:boxSDP}).
This method is an iterative method and it searches a point which satisfies 
a first-order optimality condition.

We say that $\X^* \in \FC$ satisfies a first-order optimality condition of (\ref{eq:boxSDP})
if 
\begin{eqnarray}
\langle \nabla f(\X^*) \ | \ \X - \X^* \rangle \ge 0 \qquad \mbox{for} \quad \forall \X \in \FC.
\label{eq:first-order}
\end{eqnarray}
Here,  we use $\langle \A \ | \ \B \rangle $ to denote the inner-product between  
 $\A \in \SMAT^n$ and $\B \in \SMAT^n$, 
and $\nabla f(\X^*) \in \SMAT^n$ is the gradient matrix of $f$ at $\X^*$.
In particular, when $f(\X)$ is a convex function, a point $\X^* \in \FC$ that satisfies 
(\ref{eq:first-order}) is an optimal solution.
We can derive an equivalent but more convenient condition for $\X^* \in \FC$,
\begin{eqnarray*}
 \underline{f}(\X^*) = 0
\end{eqnarray*}
where
\begin{eqnarray}
\underline{f}(\widehat{\X}) := \min \ \{ \langle \nabla f(\widehat{\X}) \ | \ \X - \widehat{\X} \rangle 
\ : \ \X \in \FC \}.
\label{eq:underlinef}
\end{eqnarray}

Xu {\it et al} \cite{XU11} proved that the feasible direction method generates an sequence
$\{\X^k\} \subset \FC$ that 
attains $\lim_{k \to \infty} \underline{f}(\X^k) = 0$.
They conducted numerical tests on 
simple objective functions that involved the variable matrix $\X$ in linear or quadratic terms.

In this paper, we propose an iterative method
for the box-constrained SDP~(\ref{eq:boxSDP})
using the distance from the current point to the boundary of $\FC$.
We introduce a concept of  the distance to the boundary of the feasible set 
from a trust-region method of Coleman and Li~\cite{COLEMAN96}
proposed for the simple-bound problem~(\ref{eq:boxLP}).
However, we can not directly apply the search direction of \cite{COLEMAN96}  to 
the box-constrained SDP~(\ref{eq:boxSDP}) 
by copying the interval condition $\l \le \x \le \u$ to 
the eigenvalue conditions $\O \preceq \X \preceq \I$,
since the matrix $\X$ 
involves not only the eigenvalues
but also the eigenvectors.
In particular, 
it is not straightforward to guarantee  a non-zero step length 
if we define a search direction ignoring the property 
that the eigenvectors are not always continuous functions on $\X$.
We devise a new search direction by taking both the eigenvalues
and the eigenvectors into consideration.
We give a non-zero range of the step length, and 
we ensure that a movement along the search direction in this range remains in $\FC$.

We also introduce a quadratic approximation function and a radius adjustment
from the trust-region methods~\cite{CONN00, GRIVA09, NOCEDAL06, TANG2014, ZHOU13}.
In ordinary trust-region methods, the search direction is obtained by solving 
a trust-region sub-problem, and the sub-problem is usually an optimization problem 
that minimizes a quadratic function 
with a constraint where the search direction is bounded by a trust-region radius.
The search direction by such a trust-region sub-problem
was examined for nonlinear semidefinite complementarity programs in \cite{LEIBFRITZ02}, but
an evaluation of the second derivative functions required a huge computation cost
and the problem size there was at most $n = 100$.
In our approach, we first obtain the search direction 
based on  the distance to the boundary, 
then we obtain the step length along this search direction so that 
the next point will stay in the region determined by a radius.
In the computation of the step length, we utilize the second derivative in its quadratic form,
hence the computation cost in each iteration is lower than the evaluation of whole entries of 
the second derivative.
We update the radius for the next iteration using an deviation of the
quadratic approximation function from the objective function.

In this paper, we discuss convergence properties of the generated sequence
for the first-order optimality condition.
Numerical tests in this paper show
that the proposed method
solves strongly-nonlinear  functions faster than the feasible direction method.
The computation cost of the proposed method in each iteration 
is low compared to the penalty barrier method implemented in
PENLAB~\cite{PENLAB13},
and the proposed method can handle larger problems than the penalty barrier method.
This paper is organized as follows.
Section~\ref{sec:trust} discusses equivalent conditions of the first-optimality conditions.
We introduce the new search direction $\D(\X)$, and propose 
the iterative method with adaptive radius adjustment
in Algorithm~\ref{al:trust}.
Section~\ref{sec:convergence} establishes the convergence properties of the proposed method.
Section~\ref{sec:results} reports numerical results on the performance comparison
of the proposed method, the feasible direction, 
and the penalty barrier method.
Finally, Section~\ref{sec:conclusions} gives a conclusion of this paper and discusses future directions.

\subsection{Notation and preliminaries} \label{sec:notation}
The inner-product between $\A \in \Real^{m \times n}$
and $\B \in \Real^{m \times n} $ 
is defined by $\langle \A | \B \rangle := Trace(\A^T \B)$.
Here, $Trace(\X)$ for a matrix $\X \in \Real^{n \times n}$ 
is the summation of its diagonal elements, that is,
$Trace(\X) := \sum_{i=1}^n X_{ii}$.

For $\A \in \Real^{m\times n}$, we define the Frobenius norm 
by $||\A||_F := \sqrt{\langle \A \ | \ \A \rangle}$.
From the Cauchy-Schwartz inequality, it holds 
$\left| \langle \A \ | \ \B \rangle  \right| \le ||\A||_F ||\B||_F$ 
for $\A \in \Real^{m \times n}$
and $\B \in \Real^{m \times n}$.
Throughout the paper, 
we often use the relation $ \langle \A \ | \ \B \rangle = \langle \B \ | \ \A \rangle$.
In addition, we use the inequality 
$\langle \A \ | \  \B \rangle \ge 0$ for two positive semidefinite matrices 
$\A \succeq \O$ and $\B \succeq \O$.

For a symmetric matrix $\A \in \SMAT^n$,  the 2-norm $||\A||_2$  is defined 
by the largest absolute eigenvalue of $\A$.
The notation $diag(\kappa_1, \kappa_2, \ldots, \kappa_n)$ stands for 
the diagonal matrix whose diagonal elements are $\kappa_1,\kappa_2, \ldots, \kappa_n$.
When $\A = \Q \K \Q^T $ is the eigenvalue decomposition of $\A$ with the diagonal matrix
$\K = diag(\kappa_1, \kappa_2, \ldots, \kappa_m)$, the $r$th power of $\A$ 
for $r \in \Real$
is given by 
$\A^r := \Q diag(\kappa_1^r, \kappa_2^r, \ldots, \kappa_n^r) \Q^T$.



The gradient matrix $\nabla f(\X) \in \SMAT^n$ and 
the Hessian mapping $\nabla^2 f(\X)$ at $\X \in \SMAT^n$ are defined so that 
a Taylor expansion for $\D \in \SMAT^n$ holds with
\begin{eqnarray*}
f(\X + \D) = f(\X) + \langle \nabla f(\X) \ | \ \D \rangle
+ \frac{1}{2} \langle \D \ | \ \nabla^2 f(\X) \ |  \ \D \rangle
+ O(||\D||_F^2),
\end{eqnarray*}
where $O(d)$ is of the order of $d$.
For example, for a function $\hat{f}(\X) = \langle \X \ | \ \X \rangle$,
we have $\nabla \hat{f}(\X) = 2 \X$ 
and $\langle \D \ | \ \nabla^2 \hat{f}(\X) \ | \ \D \rangle = 2 \langle \D \ | \ \D \rangle $
from the relation
$\langle \X + \D  \ | \ \X + \D \rangle = 
\langle \X   \ | \  \X \rangle + 2 \ \langle \X  \ | \  \D \rangle + \langle \D  \ | \  \D \rangle$.
The gradient matrix $\nabla f(\X)$ corresponds to the Fr\'{e}chet derivative,
and we have
$\langle \A \ | \ \nabla^2 f(\X) \ | \ \B \rangle = \sum_{i,j,k,l=1}^n 
\frac{\partial^2 f(\X)}{\partial X_{kl} \partial X_{ij}} A_{ij} B_{kl}$
for $\A, \B \in \SMAT^n$.

We use the matrices $\P(\X)$ and $\bold{\Gamma}(\X)$ 
to denote the eigenvalue decomposition of $\nabla f(\X)$ as
$\nabla f(\X) = \P(\X) \bold{\Gamma}(\X) \P(\X)^T$. The matrix $\bold{\Gamma}(\X)$ is 
the diagonal matrix whose diagonal elements are the descending-order eigenvalues of $\nabla f(\X)$, 
denoted by $\gamma_1(\X) \ge \gamma_2(\X) \ge \ldots \ge \gamma_n(\X)$.
The $j$th column of $\P(\X)$, denoted by $\p_j(\X)$, 
is the associated eigenvector of $\gamma_j(\X)$.
We use $n_+(\X)$ and $n_-(\X)$ to 
denote the number of positive and non-positive eigenvalues of $\nabla f(\X)$, 
respectively. We divide $\bold{\Gamma}(\X)$ into the two blocks,
$\bold{\Gamma}_+(\X) := diag(\gamma_1(\X), \gamma_2(\X), \ldots, \gamma_{n_+(\X)})$, 
$\bold{\Gamma}_-(\X) := diag(\gamma_{n_+(\X) + 1} (\X), \gamma_{n_+(\X) + 2} (\X), \ldots, \gamma_{n})$.
Note that the sizes of $\bold{\Gamma}_+(\X)$ and $\bold{\Gamma}_-(\X)$
can be zero, but the total is always $n_+(\X) + n_-(\X) = n$.
We also divide $\P(\X)$ into the two matrices $\P_+(\X), \P_-(\X)$ by collecting the 
corresponding vectors, so the columns of $\P_+(\X)$ are $\p_1(\X),\ldots,\p_{n_+(\X)}(\X)$ in this order.
As a property of eigenvectors,
we have $\P_+(\X)^T \P_-(\X) = \O$. We also know that $\P_+(\X)^T \P_+(\X)$ 
is the identity matrix of dimension $n_+(\X)$ and
$\P_-(\X)^T \P_-(\X)$ is the identity matrix of dimension $n_-(\X)$.
Finally, we define $\gamma_{\max}(\X) := || \nabla f(\X)||_2$. 
From the definition of the  2-norm, it holds that
$\gamma_{\max}(\X)  = \max\{|\gamma_1(\X)|, |\gamma_n(\X)|\}$.

\section{An iterative method using boundary distance information}\label{sec:trust}

For the simple bound problem (\ref{eq:boxLP}), 
Coleman and Li~\cite{COLEMAN96} proposed a trust-region method which
measures the distance 
from the current feasible point $\x \in \Real^n$ to the boundary of the feasible set $(\l \le \x \le \u)$. 
They defined the vector $\v(\x) \in \Real^n$ as
\begin{eqnarray*}
v_i(\x) = \left\{
\begin{array}{lll}
x_i - l_i & \mbox{if} & \frac{\partial f(\x)}{\partial x_i} > 0 \\
u_i - x_i & \mbox{if} & \frac{\partial f(\x)}{\partial x_i} \le 0.
\end{array}\right.
\end{eqnarray*}
This vector was used to control the approach to the boundary,
and the key observation in the discussion of \cite{COLEMAN96}
was that  $\x^*$ satisfies the first-order optimality condition if and only if
$\frac{\partial f(\x)}{\partial x_i} v_i(\x) = 0$ for each $i = 1,\ldots, n$.

We can not directly extend the definition of $\v(\x)$ to
the box-constrained SDPs (\ref{eq:boxSDP}) using the conditions on the eigenvalue of $\X$, since 
the distance to the boundary of $\FC$
relates to not only the eigenvalues but also the eigenvectors.
To take the effect of eigenvectors into account, 
we define two positive semidefinite matrices for $\X \in \FC$;
\begin{eqnarray*}
\V_+(\X)  := \P_+(\X)^T \X \P_+(\X), \quad \mbox{and} \quad 
\V_-(\X)  := \P_-(\X)^T (\I - \X) \P_-(\X).
\end{eqnarray*}
The definition of these matrices brings us other properties of
the first-order optimality condition
in Lemma~\ref{le:equivalent}.
In the lemma,
we use a matrix $\D(\X) \in \SMAT^n$ 
and a scalar $N(\X)$ defined by
\begin{eqnarray}
\D(\X) &:=& \P(\X) 
\left(\begin{array}{cc}
\V_+(\X)^{1/2} \bold{\Gamma}_+(\X) \V_+(\X)^{1/2} & \gamma_{\max}(\X) \P_+(\X)^T \X \P_-(\X) \\
\gamma_{\max}(\X) \P_-(\X)^T \X \P_+(\X) & \V_-(\X)^{1/2} \bold{\Gamma}_-(\X) \V_-(\X)^{1/2} 
\end{array}\right)
\P(\X)^T \label{eq:D}\\
N(\X) &:=& \langle \nabla f(\X) \ | \ \D(\X) \rangle. \label{eq:N}
\end{eqnarray}
The definition of the matrix $\D(\X)$ includes
the distance information to the boundary of the feasible sets $\FC$ 
via the matrices $\V_+(\X)$ and $\V_-(\X)$ like $\v(\x)$ above.

Using the relations
$\nabla f(\X)  
= \P_+(\X) \bold{\Gamma}_+(\X) \P_+(\X)^T +  \P_-(\X) \bold{\Gamma}_-(\X) \P_-(\X)^T$,
we can compute $||\D(\X)||_F^2$ and $N(\X)$ as follow;
\begin{eqnarray}
||\D(\X)||_F^2 &=& 
||\V_+(\X)^{1/2} \bold{\Gamma}_+(\X) \V_+(\X)^{1/2}||_F^2 
+ ||\V_-(\X)^{1/2} \bold{\Gamma}_-(\X) \V_-(\X)^{1/2}||_F^2 \nonumber \\
& & + 2 \gamma_{\max}^2 ||\P_{+}(\X)^T \X \P_-(\X)||_F^2,  \label{eq:anotherD} \\
N(\X) &=& ||\V_+(\X)^{1/4} \bold{\Gamma}_+(\X) \V_+(\X)^{1/4}||_F^2 
+ ||\V_-(\X)^{1/4} \bold{\Gamma}_-(\X) \V_-(\X)^{1/4}||_F^2. \label{eq:anotherN}
\end{eqnarray}

\begin{LEMM}\label{le:equivalent}

For a matrix $\X^* \in \FC$, the following conditions are equivalent.
\begin{enumerate}[(a)]
\item $\X^*$ satisfies the first-order optimality condition (\ref{eq:first-order}).
\item  $\langle \bold{\Gamma}_+(\X^*) \ | \  \V_+(\X^*) \rangle 
= \langle \bold{\Gamma}_-(\X^*) \ | \ \V_-(\X^*) \rangle = 0.$ 
\item $N(\X^*) = 0.$
\item $||\D(\X^*)||_F = 0.$
\end{enumerate}
\end{LEMM} 

\noindent{\bf Proof:}

\noindent[$(a) \Rightarrow (b)$] 
We define a matrix $\widehat{\X} := \P_+(\X^*) \P_+(\X^*)^T \X^* \P_+(\X^*) \P_+(\X^*)^T 
+ \P_-(\X^*)  \P_-(\X^*)^T$.
Since $\X^* \in \FC$, we obtain $\widehat{\X} \succeq \O$ and 
\begin{eqnarray*}
\I - \widehat{\X} &=& (\P_+(\X^*) \P_+(\X^*)^T + \P_-(\X^*) \P_-(\X^*)^T) \\
& &  - 
(\P_+(\X^*) \P_+(\X^*)^T \X^* \P_+(\X^*) \P_+(\X^*)^T 
 + \P_-(\X^*)  \P_-(\X^*)^T) \\
&=& \P_+(\X^*) \P_+(\X^*)^T (\I-  \X^* ) \P_+(\X^*) \P_+(\X^*)^T \succeq \O,
\end{eqnarray*}
hence, $\widehat{\X} \in \FC$.
Substituting $\widehat{\X} \in \FC$ into the inequality (\ref{eq:first-order}),
we have
\begin{eqnarray*}
& & \langle \nabla f(\X^*)  \ | \ \widehat{\X}  - \X^* \rangle  \\
&=& \langle \P_+(\X^*) \bold{\Gamma}_+(\X^*) \P_+(\X^*)^T + \P_-(\X^*) \bold{\Gamma}_-(\X^*) \P_-(\X^*) ^T \\
& & \qquad
\ | \ \P_+(\X^*) \P_+(\X^*)^T \X^* \P_+(\X^*) \P_+(\X^*)^T + \P_-(\X^*)  \P_-(\X^*)^T - \X^* \rangle  \\
&=& \langle \bold{\Gamma}_-(\X^*) \ | \ \I  \rangle
- \langle \bold{\Gamma}_-(\X^*) \ | \ \P_-(\X^*)^T \X^* \P_-(\X^*) \rangle 
= \langle \bold{\Gamma}_-(\X^*) \ | \ \V_-(\X^*)  \rangle \ge 0.
\end{eqnarray*}
Here, we used
$\langle \A \ | \ \B \rangle = Trace(\A^T\B) = Trace(\B^T \A)$, 
$\P_+(\X^*)^T \P_+(\X^*) = \I$ \newline 
and $\P_+(\X^*)^T \P_-(\X^*) = \O$.
Since $- \bold{\Gamma}_-(\X^*) \succeq \O$ and 
$\V_-(\X^*) \succeq \O$, we also have \newline
$\langle -\bold{\Gamma}_-(\X^*) \ | \ \V_-(\X^*)  \rangle \ge 0$, so that
we obtain $\langle  \bold{\Gamma}_-(\X^*) \ | \V_-(\X^*) \rangle = 0$.

Similarly, for the matrix
$\overline{\X} :=  \P_-(\X^*)  \P_-(\X^*)^T \X^* \P_-(\X^*)  \P_-(\X^*)^T \succeq \O$,
we can show $\I - \overline{\X} = \P_+(\X^*) \P_+(\X^*)^T + \P_-(\X^*) \P_-(\X^*)^T (\I - \X^*) \P_-(\X^*) \P_-(\X^*)^T \succeq \O$, therefore we have $\overline{\X} \in \FC$.
Putting $\overline{\X} \in \FC$ into  (\ref{eq:first-order}), we have
$\langle \nabla f(\X^*)  \ | \ \overline{\X}  - \X^* \rangle = 
- \langle \bold{\Gamma}_+(\X^*) \ | \ \V_+(\X^*)  \rangle \ge 0$. 
On the other hand, from the properties
$\bold{\Gamma}_+(\X^*) \succeq \O$ and  $\V_+(\X^*) \succeq \O$, 
it holds $\langle  \bold{\Gamma}_+(\X^*) \ | \V_+(\X^*) \rangle \ge 0$.
Hence, we obtain $\langle  \bold{\Gamma}_+(\X^*) \ | \V_+(\X^*) \rangle = 0$. 

\noindent[$(b) \Rightarrow (a)$] 
For any $\X \in \FC$,
it holds that 
\begin{eqnarray*}
 & & \langle \nabla f(\X^*) \ |\ \X - \X^* \rangle \\
&=& \langle 
\P_+(\X^*) \bold{\Gamma}_+(\X^*)\P_+(\X^*)^T 
+ \P_-(\X^*) \bold{\Gamma}_-(\X^*) \P_-(\X^*)^T  \ | \ \X - \X^* \rangle \\
&=& \langle \bold{\Gamma}_+(\X^*)  \ | \ \P_+(\X^*)^T \X  \P_+(\X^*) \rangle
- \langle  \bold{\Gamma}_+(\X^*)  \ | \ \V_+(\X^*) \rangle   \\
& & 
\qquad - \langle  \bold{\Gamma}_-(\X^*)  \ |  \ \P_-(\X^*)^T (\I-\X)  \P_-(\X^*) \rangle 
+ \langle  \bold{\Gamma}_-(\X^*)  \ | \ \V_-(\X^*) \rangle  \\
&=& \langle \bold{\Gamma}_+(\X^*)  \ | \  \P_+(\X^*)^T \X  \P_+(\X^*) \rangle
+ \langle - \bold{\Gamma}_-(\X^*)  \ |  \ \P_-(\X^*)^T (\I-\X)  \P_-(\X^*) \rangle \ \ge 0.
\end{eqnarray*}
For the last equality, we used $\langle  \bold{\Gamma}_+(\X^*)  \ | \ \V_+(\X^*) \rangle  = 
\langle  \bold{\Gamma}_-(\X^*) \ | \ \V_-(\X^*) \rangle  = 0 $ from $(b)$.
In addition, the last non-negativity came from
$\P_+(\X^*)^T \X  \P_+(\X^*) \succeq \O$ and $\P_-(\X^*)^T (\I-\X)  \P_-(\X^*) \succeq \O$.

\noindent[$(b) \Rightarrow (c)$] 
Since $\langle \bold{\Gamma}_+(\X^*) \ | \ \V_+(\X^*) \rangle 
= Trace(\V_+(\X^*)^{1/2} \bold{\Gamma}_+(\X^*) \V_+(\X^*)^{1/2}) $
and \newline $\V_+(\X^*)^{1/2} \bold{\Gamma}_+(\X^*) \V_+(\X^*)^{1/2} \succeq \O$,
the condition
$\langle \bold{\Gamma}_+(\X^*) \ | \ \V_+(\X^*) \rangle = 0$ indicates
all the eigenvalues of 
$\V_+(\X^*)^{1/2} \bold{\Gamma}_+(\X^*) \V_+(\X^*)^{1/2}$ are 0, therefore, 
$\V_+(\X^*)^{1/2} \bold{\Gamma}_+(\X^*) \V_+(\X^*)^{1/2} 
= \O$. 
We now consider the eigenvalue decomposition $\V_+(\X^*) = \Q \K \Q^T$ 
such that \newline
$\K = diag(\kappa_1, \kappa_2, \ldots, \kappa_{n_+(\X^*)})$ is 
the diagonal matrix with the eigenvalues of $\V_+(\X^*)$.
Since $\V_+(\X^*) \succeq \O$, it holds that 
$\kappa_i  \ge 0$ for $i = 1, \ldots, n_+(\X^*)$.
We define a positive semidefinite matrix $\W := \Q^T \bold{\Gamma}_+(\X^*) \Q$.
Since the matrix $\Q$ is an orthogonal matrix, 
$\V_+(\X^*)^{1/2} \bold{\Gamma}_+(\X^*) \V_+(\X^*)^{1/2} = \O$
leads to $\K^{1/2} \W \K^{1/2} = \O$.
By taking the diagonal elements, we know $\kappa_i^{1/2}W_{ii}\kappa_i^{1/2} = 0$
for $i = 1, \ldots, n_+(\X^*)$.
Therefore, it holds that $\kappa_i^{1/4}W_{ii}\kappa_i^{1/4} = 0$.
Since a matrix $\K^{1/4} \W \K^{1/4}$
is positive semidefinite and its diagonal elements are zero,
we obtain $\K^{1/4} \W \K^{1/4} = \O$, hence,
$\V_+(\X^*)^{1/4} \bold{\Gamma}_+(\X^*) \V_+(\X^*)^{1/4} = \O$.
Similarly,
the condition 
$\langle \bold{\Gamma}_-(\X^*) \ | \ \V_-(\X^*) \rangle = 0$ implies
 $\V_-(\X^*)^{1/4} \bold{\Gamma}_-(\X^*) \V_-(\X^*)^{1/4} = \O$.
Hence, 
we obtain $(c)$ by (\ref{eq:anotherN}).

\ifx10
より詳しくは、$\V_+^{1/2} = \Q\M\Q^{T}$ と固有値分解すると、$\A = \Q^T \bold{\Gamma}_+ \Q$ で
\begin{eqnarray*}
& & \V_+^{1/2} \bold{\Gamma}_+ \V_+^{1/2}  = \O \\
&\Leftrightarrow& \Q\M\Q^T \bold{\Gamma}_+ \Q\M\Q^T = \O \\
&\Leftrightarrow& \M\A\M = \O \\
&\Leftrightarrow& \mu_i\mu_j A_{ij} = 0 \\
&\Leftrightarrow& \sqrt{\mu_i}\sqrt{\mu_j}A_{ij} = 0 \\
&\Leftrightarrow&  \M^{1/2}\A\M^{1/2} = \O \\
&\Leftrightarrow & \V_+^{1/4} \bold{\Gamma}_+ \V_+^{1/4}  = \O
\end{eqnarray*}
\fi

\noindent[$(c) \Rightarrow (b)$]
The condition $N(\X^*) = 0$ leads to 
$\V_+(\X^*)^{1/4} \bold{\Gamma}_+(\X^*) \V_+(\X^*)^{1/4} = \O$
\newline
and $\V_-(\X^*)^{1/4} \bold{\Gamma}_-(\X^*) \V_-(\X^*)^{1/4} = \O$.
Hence, it holds that
\begin{eqnarray*}
\langle \bold{\Gamma}_+(\X^*) \ | \ \V_+(\X^*) \rangle &=&
\langle \bold{\Gamma}_+(\X^*) \ | \ \V_+(\X^*)^{1/4} \V_+(\X^*)^{1/2} \V_+(\X^*)^{1/4}  \rangle \\
&=& \langle \V_+(\X^*)^{1/4} \bold{\Gamma}_+(\X^*) \V_+(\X^*)^{1/4} 
\ | \  \V_+(\X^*)^{1/2}  \rangle = 0.
\end{eqnarray*}
Similarly, we obtain 
$\langle \bold{\Gamma}_-(\X^*) \ | \ \V_-(\X^*) \rangle = 0$ from
$\V_-(\X^*)^{1/4} \bold{\Gamma}_-(\X^*) \V_-(\X^*)^{1/4} = \O$.

\noindent[$(b) \Rightarrow (d)$] 
As a first step of [$(b) \Rightarrow (c)$] above,
we obtained $\V_+(\X^*)^{1/2} \bold{\Gamma}_+(\X^*) \V_+(\X^*)^{1/2}  = \O$  and \newline
$\V_-(\X^*)^{1/2} \bold{\Gamma}_-(\X^*) \V_-(\X^*)^{1/2}  = \O$.
Since all the eigenvalues in $\bold{\Gamma}_+(\X^*)$ are positive,
the properties $\langle \bold{\Gamma}_+(\X^*) \ | \ \V_+(\X^*) \rangle  = 0$ and $\V_+(\X^*) \succeq \O$
lead to
$\V_+(\X^*) = \O$. Furthermore, the decomposition $\V_+(\X^*) = \P_+(\X^*)^T (\X^*)^{1/2} (\X^*)^{1/2} \P_+(\X^*) = \O$
implies \newline
$ \P_+(\X^*)^T (\X^*)^{1/2}  = \O$.
Therefore, it holds that \newline
$\P_+(\X^*)^T \X^* \P_-(\X^*) = \P_+(\X^*)^T (\X^*)^{1/2} (\X^*)^{1/2} \P_-(\X^*) = \O$.
Hence, we conclude $||\D(\X^*)||_F = 0$ from (\ref{eq:anotherD}).

\noindent[$(d) \Rightarrow (b)$] 
From the relation (\ref{eq:anotherD}), the condition $||\D(\X^*)||_F = 0$ indicates \newline
$\V_+(\X^*)^{1/2} \bold{\Gamma}_+(\X^*) \V_+(\X^*)^{1/2}  = \O$ and $\V_-(\X^*)^{1/2} \bold{\Gamma}_-(\X^*) \V_-(\X^*)^{1/2}  = \O$.
By taking the traces of these matrices, we obtain $(b)$.
\hfill$\square$

Lemma~\ref{le:equivalent}, (\ref{eq:N}) and (\ref{eq:anotherN}) indicate that
when $\X$ does not satisfy the first-order optimality condition,
we can take $- \frac{\D(\X)}{||\D(\X)||_F} $ as a descent direction of $f(\X)$, 
that is, $ \langle \nabla f(\X) \ | \ -\frac{\D(\X)}{||\D(\X)||_F} \rangle < 0$.
Hence, we can expect that the decrease of the objective function
$f(\X - \alpha \frac{\D(\X)}{||\D(\X)||_F})  < f(\X)$
for a certain value $\alpha > 0$.
The next lemma gives a non-zero range of $\alpha$ to ensure
$\X - \alpha \frac{\D(\X)}{||\D(\X)||_F}  \in \FC$.

\begin{LEMM}\label{le:rangealpha}
If $\X \in \FC$ does not satisfy the first-order optimality condition, then
$\X - \alpha \frac{\D(\X)}{||\D(\X)||_F} \in \FC$ for 
$\alpha \in [0,\frac{||\D(\X)||_F}{\gamma_{\max}(\X)}]$.
\end{LEMM}

\noindent{\bf Proof:}

From the definition of $\gamma_{\max}(\X)$, 
the matrix $\I - \frac{\bold{\Gamma}_+(\X)}{\gamma_{\max}(\X)}$ is a nonnegative diagonal matrix, 
hence this matrix is positive semidefinite.
Using $\P(\X)\P(\X)^T = \I$ and $-\bold{\Gamma}_-(\X) \succeq \O$,
it holds
{\footnotesize
\begin{eqnarray*}
& & \X - \frac{\D(\X)}{\gamma_{\max}(\X)} 
= \P(\X) \P(\X)^T \left( \X - \frac{\D(\X)}{\gamma_{\max}(\X)}  \right) \P(\X) \P(\X)^T \\
&=& \P(\X) \left\{
\P(\X)^T \X \P(\X) - 
\left(\begin{array}{cc}
\V_+(\X)^{1/2}\frac{\bold{\Gamma}_+(\X)}{\gamma_{\max}(\X)}\V_+(\X)^{1/2} & \P_+(\X)^T \X \P_-(\X) \\
\P_-(\X)^T \X \P_+(\X) & \V_-(\X)^{1/2}\frac{\bold{\Gamma}_-(\X)}{\gamma_{\max}(\X)}\V_-(\X)^{1/2} 
\end{array}\right) \right\} \P(\X)^T \\
&=& \P(\X)
\left(\begin{array}{cc}
\V_+(\X)^{1/2}\left(\I - \frac{\bold{\Gamma}_+(\X)}{\gamma_{\max}(\X)}\right) \V_+(\X)^{1/2} & \O \\
\O & \P_-(\X)^T \X \P_-(\X) + \V_-(\X)^{1/2}  \frac{(-\bold{\Gamma}_-(\X))}{\gamma_{\max}(\X)}\V_-(\X)^{1/2} 
\end{array}\right)  \P(\X)^T \\
& \succeq  & \O.
\end{eqnarray*}
}

In a similar way, noticing $\P_+(\X)^T (\I-\X) \P_-(\X) = - \P_+(\X)^T \X \P_-(\X)$ 
and $\I + \frac{\bold{\Gamma}_-(\X)}{\gamma_{\max}(\X)} \succeq \O$, we derive
{\footnotesize
\begin{eqnarray*}
& & \I - \X + \frac{\D(\X)}{\gamma_{\max}(\X)} \\
&=& \P(\X)
\left(\begin{array}{cc}
\P_+(\X)^T (\I-\X) \P_+(\X) + \V_+(\X)^{1/2}  \frac{\bold{\Gamma}_+(\X)}{\gamma_{\max}(\X)}\V_+(\X)^{1/2}  & \O \\
\O & \V_-(\X)^{1/2}\left(\I + \frac{\bold{\Gamma}_-(\X)}{\gamma_{\max}(\X)}\right) \V_-(\X)^{1/2}
\end{array}\right)  \P(\X)^T \\
& \succeq  & \O.
\end{eqnarray*}
}

From two linear combinations
\begin{eqnarray*}
\X - \alpha \frac{\D(\X)}{||\D(\X)||_F} &=& 
\left( 1 - \alpha \frac{\gamma_{\max}(\X)}{||\D(\X)||_F}\right) \X 
+ \alpha \frac{\gamma_{\max}(\X)}{||\D(\X)||_F}
\left( \X - \frac{\D(\X)}{\gamma_{\max}(\X)} \right) \\
\I - \left(\X - \alpha \frac{\D(\X)}{||\D(\X)||_F} \right) &=&
\left( 1 - \alpha \frac{\gamma_{\max}(\X)}{||\D(\X)||_F}\right) (\I-\X)
+ \alpha \frac{\gamma_{\max}(\X)}{||\D(\X)||_F}
\left( \I - \X + \frac{\D(\X)}{\gamma_{\max}(\X)} \right),
\end{eqnarray*}
we obtain $\X - \alpha \frac{\D(\X)}{||\D(\X)||_F} \succeq \O$ and
$\I - \left(\X - \alpha \frac{\D(\X)}{||\D(\X)||_F}\right) \succeq \O$
for $\alpha \in [0,\frac{||\D(\X)||_F}{\gamma_{\max}(\X)}]$.

\hfill$\square$

Based on the property that $-\frac{\D(\X)}{||\D(\X)||_F}$ is a descent direction of $f(\X)$, 
we can use $\S(\X) := \frac{\D(\X)}{||\D(\X)||_F}$
as a normalized search direction to find a minimizer.

We propose an iterative method for the box-constrained SDP (\ref{eq:boxSDP}) 
as Algorithm~\ref{al:trust}.
In Algorithm~\ref{al:trust}, we use
a quadratic approximation of $f$ with the direction $\S(X)$;
\begin{eqnarray*}
q(\alpha,\X) := 
f(\X) - \alpha \langle \nabla f(\X) \ | \ \S(\X) \rangle
+ \frac{\alpha^2}{2} \langle \S(\X) \ | \ \nabla^2 f(\X) \ | \ \S(\X) \rangle.
\end{eqnarray*}
We should note that the quadratic approximation function $q(\alpha, \X^k)$
requires $\nabla^2 f(\X^k)$ in only the scalar value $\langle \S(\X^k) \ | \ \nabla^2 f(\X^k) \ | \ \S(\X^k) \rangle$.
Hence, we do not always need to evaluate each element of $\nabla^2 f(\X^k)$
 in each iteration.
For example, for a function $\hat{f}(\X) = \cos(\langle \X \ | \ \X \rangle)$ and
a symmetric matrix $\S \in \SMAT^n$, 
it holds $\langle \S \ | \ \nabla^2 \hat{f}(\X) \ | \ \S \rangle
 = -2\sin(\langle \X \ | \ \X \rangle) \langle \S \ | \ \S \rangle
-4\cos(\langle \X \ | \ \X \rangle) \langle \X \ | \ \S \rangle^2$.
This makes each iteration of  Algorithm~\ref{al:trust} low cost compared to
the interior-point methods or the penalty barrier method.

We note that the generated sequence by Algorithm~\ref{al:trust}
remains in $\FC$, that is, $\{ \X^k \} \subset \FC$ from Lemma~\ref{le:rangealpha}.
In Steps 4 and 5, we adjust the radius $\Delta_k$.
This adjustment is necessary to discuss the convergence properties.

\begin{table}[t]
\noindent \hrulefill
\begin{ALGO}\label{al:trust} {\rm An iterative method using boundary distance  
for box-constrained SDPs}
{\rm
\begin{enumerate}[Step 1:]
\item Choose an initial point $\X^0 \in \FC$. 
Set an initial radius $\Delta_0 >0$ and set a stopping threshold $\epsilon > 0$.
Choose parameters $\mu_1, \mu_2, \eta_1, \eta_2$ such that
$0 < \mu_1 < \mu_2$ and $0 < \eta_1 < 1 < \eta_2$. 
Set an iteration count $k=0$. 
\item If $N(\X^k) < \epsilon$, output $\X^k$ as a solution and stop.
\item Solve a quadratic problem with respect to $\alpha$;
\begin{eqnarray}
\min\quad 
q(\alpha, \X^k)
\quad \mbox{subject to} \quad 0 \le \alpha \le 
\min\left\{ \frac{||\D(\X^k)||_F}{\gamma_{\max}(\X^k)}, \Delta_k \right\},
\label{eq:subproblem}
\end{eqnarray}
and let the step length $\alpha_k$ be the minimizer of (\ref{eq:subproblem}).
\item Let $\overline{\X}^k := \X^k - \alpha_k \S(\X^k)$ where
$\S(\X^k) := \frac{\D(\X^k)}{||\D(\X^k)||_F}$. 
Compute the ratio 
\begin{eqnarray}
r_k := \frac{f(\X^k) - f(\overline{\X}^k)}{f(\X^k) - q(\alpha_k, \X^k)},
\end{eqnarray}
and set 
\begin{eqnarray*}
\X^{k+1} = \left\{
\begin{array}{ll}
\overline{\X}^k & \mbox{if} \quad
r_k \ge \mu_1\\
\X^k & \mbox{otherwise}.
\end{array}\right.
\end{eqnarray*}
\item Update the radius $\Delta_{k}$ by 
\begin{eqnarray*}
\Delta_{k+1} = \left\{
\begin{array}{ll}
\eta_1 \Delta_k & \mbox{if} \quad r_k < \mu_1 \\
 \multicolumn{1}{r}{\Delta_k} & \mbox{if} \quad \mu_1 \le r_k \le \mu_2 \\
\eta_2 \Delta_k & \mbox{if} \quad r_k > \mu_2. \\
\end{array}\right.
\end{eqnarray*}
\item Set $k \leftarrow k+1$ and return to Step 2.
\end{enumerate}
}
\end{ALGO}
\noindent \hrulefill
\end{table}


\section{Convergence properties}\label{sec:convergence}

A matrix $\X^* \in \FC$ satisfies 
the first-order optimality condition (\ref{eq:first-order}) if and only if 
$\underline{f}(\X^*) = 0$, as noted in Section~1.
In this section, we show that 
the sequence $\{\X^k\} \subset \FC$ generated by Algorithm~\ref{al:trust} 
with the stopping threshold $\epsilon = 0$
attains $\lim_{k \to \infty} \underline{f}(\X^k) = 0$.
We divide the proof into two parts.
The first part shows there exists 
a subsequence of $\{N(\X^k)\}$ that converges to zero. 
The second part shows $\lim_{k \to \infty} N(\X^k) = 0$ in Theorem~\ref{th:whole},
and finally 
$\lim_{k\to \infty} \underline{f}(\X_k) = 0$ in Theorem~\ref{th:flimit}.

Using the matrix $\D(\X)$, we can employ similar approaches to \cite{COLEMAN96} 
for the proof of the first part.
However, we can not directly apply the results of \cite{COLEMAN96} to the second part. 
This is mainly because
that the eigenvector matrices $\P_+(\X)$ and $\P_-(\X)$
are not always continuous functions in $\X$. Instead, our proof relies on the 
boundedness of
$\langle \bold{\Gamma}_+(\X) \ | \ \V_+(\X) \rangle$  
and $\langle - \bold{\Gamma}_-(\X) \ | \ \V_-(\X) \rangle$.

\subsection{Convergence of subsequence}\label{sec:sub}

To analyze Algorithm~\ref{al:trust}, we introduce two constant values
\begin{eqnarray*}
M_1 &:=& \max_{\X \in \FC} ||\nabla f(\X)||_2, \\
M_2 &:=& \max_{\X \in \FC, \D \in \SMAT^n, \D \ne \O} 
\left| \frac{\langle \D \ | \ \nabla^2 f(\X) \ | \ \nabla \D \rangle}{\langle \D \ | \ \D \rangle} \right|.
\end{eqnarray*}
These values are finite from the assumptions 
that the feasible set $\FC$ is a bounded and closed set
and that the objective function $f(\X)$ 
is a twice continuously differentiable function on an open set containing $\FC$.
We can assume that $M_1 > 0$ and $M_2 > 0$ without loss of generality.
If $M_1 = 0$, then $f(\X)$ is a constant function
and every point $\X \in \FC$ is optimal. If $M_2 = 0$, then
$\nabla f(\X)$ is a constant matrix on $\FC$, so that the global minimizer can be obtained as
$\X^* = \P_-(\X) \P_-(\X)^T$ from any point $\X \in \FC$.

We now evaluate the quadratic approximation function $q(\alpha_k, \X^k)$.
\begin{LEMM}\label{le:progress}
The step length $\alpha_k$ in Step 3 satisfies 
\begin{eqnarray*}
q(\alpha_k, \X^k) \le f(\X^k)
- \frac{1}{2} \min \left\{\frac{N(\X^k)^2}{M_2 ||\D(\X^k)||_F^2},
\frac{N(\X^k)}{\gamma_{\max}(\X^k)},
\frac{\Delta_k N(\X^k)}{||\D(\X^k)||_F}  \right\}.
\end{eqnarray*}
\end{LEMM}

\noindent{\bf Proof:}



We define a quadratic function 
$\phi(\alpha) := -\alpha \frac{N(\X^k)}{||\D(\X^k)||_F}  + \frac{\alpha^2}{2} M_2$. 
From the definitions of $N(\X^k)$ and $M_2$, we have
$q(\alpha, \X^k) \le f(\X^k) + \phi(\alpha)$, hence, 
 \begin{eqnarray*}
 q(\alpha_k, \X^k) \le f(\X^k) 
 + \min_{\alpha \in \left[ 0, 
       \min\left\{ \frac{||\D(\X^k)||_F}{\gamma_{\max}(\X^k)}, \Delta_k \right\}
\right]} \phi(\alpha).
\end{eqnarray*}

Since $N(\X^k) = \langle \nabla f(\X^k) \ | \ \D(\X^k) \rangle  > 0$ 
(otherwise, 
Lemma~\ref{le:equivalent} indicates that $\X^k$ already satisfies the first-order optimality condition)
and $\phi(\alpha)$ is a quadratic function with respect to $\alpha$,
the minimum of $\phi$ is attained at one of the three candidates
$\frac{||\D(\X^k)||_F}{\gamma_{\max}(\X^k)}, \Delta_k$ or
$\widehat{\alpha} := \frac{N(\X^k)}{M_2||\D(\X^k)||_F}$.
Let $\alpha_{\min}$ be the minimizer of $\phi(\alpha)$ subject to
$0 \le \alpha \le \min\left\{ \frac{||\D(\X^k)||_F}{\gamma_{\max}(\X^k)}, \Delta_k \right\}$.


If $\alpha_{\min} = \widehat{\alpha}$,
we have 
$\phi(\widehat{\alpha}) = -\frac{1}{2} \frac{N(\X^k)^2}{M_2 ||\D(\X^k)||_F^2}$.
Next, if $\alpha_{\min} = \frac{||\D(\X^k)||_F}{\gamma_{\max}}$, we have
$\frac{||\D(\X^k)||_F}{\gamma_{\max}(\X^k)} \le \widehat{\alpha}$, therefore, 
$\frac{||\D(\X^k)||_F^2}{\gamma_{\max}(\X^k)} M_2 \le N(\X^k)$.
Hence, it holds that 
$\phi\left(\frac{||\D(\X^k)||_F}{\gamma_{\max}(\X^k)}\right) = 
-\frac{N(\X^k)}{\gamma_{\max}(\X^k)}
+ \frac{1}{2} \frac{||\D(\X^k)||_F^2}{\gamma_{max}(\X^k)^2}M_2
\le -\frac{1}{2}\frac{N(\X^k)}{\gamma_{\max}(\X^k)}$.
Finally, if $\alpha_{\min} = \Delta_k$,  the inequality $\Delta_k \le \widehat{\alpha}$ indicates that 
$\Delta_k \le \frac{N(\X^k)}{M_2||\D(\X^k)||_F}$.
Hence, it holds that 
$\phi(\Delta_k) = 
-\Delta_k \frac{N(\X^k)}{||\D(\X^k)||_F} + \frac{1}{2} \Delta_k^2 M_2
\le -\Delta_k \frac{N(\X^k)}{||\D(\X^k)||_F} + \frac{1}{2} \Delta_k \frac{N(\X^k)}{||\D(\X^k)||_F}
\le -\frac{1}{2}\frac{\Delta_k N(\X^k)}{||D(\X^k)||_F}$.

Taking the maximum of the three cases, we obtain the inequality of this lemma.
\hfill$\square$

To simplify the inequality of Lemma~\ref{le:progress},
we replace $\gamma_{\max}(\X^k)$ and $||\D(\X^k)||_F$ 
by convenient upper bounds.
Since $\gamma_{\max}(\X^k)$ is bounded by $M_1$, we consider an upper bound 
on $||\D(\X^k)||_F$. 
\begin{LEMM}\label{le:Dub}
For $\X \in \FC$, it holds that $||\D(\X)||_F^2 \le N(\X) + \frac{1}{2} M_1^2 n^3$.
\end{LEMM}

\noindent{\bf Proof:}
Let $\V_+(\X) = \Q \K \Q^T$ be the eigenvalue decomposition of $\V_+(\X)$
such that $\K = diag(\kappa_1, \kappa_2, \ldots, \kappa_{n+(\X)})$ is 
the diagonal matrix with the eigenvalues of $\V_+(\X)$.
Since $\O \preceq \X \preceq \I$, we have  $\O \preceq \V_+(\X) \preceq \I$, 
hence, $0 \le \kappa_i \le 1$ for $i = 1,2\ldots,n_+(\X)$.
Using a matrix 
$\W := \Q^T \bold{\Gamma}_+(\X) \Q$,
we compare $||\V_+(\X)^{1/2} \bold{\Gamma}_+(\X) \V_+(\X)^{1/2}||_F$
and $||\V_+(\X)^{1/4} \bold{\Gamma}_+(\X) \V_+(\X)^{1/4}||_F$;
\begin{eqnarray*}
& & ||\V_+^{1/4}(\X) \bold{\Gamma}_+(\X) \V_+^{1/4}(\X)||_F^2 - ||\V_+(\X)^{1/2} \bold{\Gamma}_+(\X) \V_+(\X)^{1/2}||_F^2 \\
&=& \langle \bold{\Gamma}_+(\X) \ |  \ \V_+(\X)^{1/2} \bold{\Gamma}_+(\X) \V_+(\X)^{1/2} \rangle \\
& & - \langle \V_+(\X)^{1/2} \bold{\Gamma}_+(\X) \V_+(\X)^{1/2} \ |  \ \V_+(\X)^{1/2} \bold{\Gamma}_+(\X) \V_+(\X)^{1/2} \rangle  \\
&=& \langle \bold{\Gamma}_+(\X) - \V_+(\X)^{1/2} \bold{\Gamma}_+(\X) \V_+(\X)^{1/2} 
\ |  \ \V_+(\X)^{1/2} \bold{\Gamma}_+(\X) \V_+(\X)^{1/2} \rangle  \\
&=& \langle \bold{\Gamma}_+(\X) - \Q \K^{1/2} \Q^T  \bold{\Gamma}_+(\X) \Q \K^{1/2} \Q^T 
\ |  \ \Q \K^{1/2} \Q^T  \bold{\Gamma}_+(\X) \Q \K^{1/2} \Q^T  \rangle  \\
&=& \langle \W - \K^{1/2} \W \K^{1/2} \ | \ \K^{1/2} \W \K^{1/2} \rangle \\
&=& ||\K^{1/4} \W \K^{1/4}||_F^2 - ||\K^{1/2} \W \K^{1/2}||_F^2  \\
&=& \sum_{i=1}^{n_+(\X)} \sum_{j=1}^{n_+(\X)} (W_{ij} \kappa_i^{1/4}\kappa_j^{1/4})^2 
- \sum_{i=1}^{n_+(\X)} \sum_{j=1}^{n_+(\X)} (W_{ij} \kappa_i^{1/2}\kappa_j^{1/2})^2  \\
&=& \sum_{i=1}^{n_+(\X)} \sum_{j=1}^{n_+(\X)} W_{ij}^2 (\kappa_i^{1/2}\kappa_j^{1/2} - \kappa_i\kappa_j)
 \ge 0.
\end{eqnarray*}
The last inequality comes from 
$0 \le \kappa_ i \le 1$ for $i = 1, \ldots, n_+(\X)$.
In a similar way, we also derive 
$||\V_-(\X)^{1/2} \bold{\Gamma}_-(\X) \V_-(\X)^{1/2}||_F^2 
\le ||\V_-(\X)^{1/4} \bold{\Gamma}_-(\X) \V_-(\X)^{1/4}||_F^2$.
We evaluate the last term of (\ref{eq:anotherD}) by a property of the Frobenius norm,
\begin{eqnarray*}
||\P_+(\X)^T \X \P_-(\X)||_F^2 \le ||\P_+(\X)||_F^2 \cdot ||\X||_F^2 \cdot ||\P_-(\X)||_F^2
\le n_+(\X) \cdot n \cdot n_-(\X) \le \frac{n^3}{4}.
\end{eqnarray*}
Here, we used
$||\P_+(\X)||_F^2 = Trace(\P_+(\X)^T \P_+(\X)) = n_+(\X)$.
In addition, we used $||\X||_F \le \sqrt{n} ||\X||_2$ from \cite[(1.2.27)]{YUAN06}
and $\O \preceq \X \preceq \I$ to derive $||\X||_F^2 \le n$,
and we used the relation 
$n_+(\X) + n_-(\X) = n$ to derive $n_+(\X) \cdot n_-(\X) \le \frac{n^2}{4}$.

Consequently, it holds from (\ref{eq:anotherD}) that 
\begin{eqnarray*}
||\D(\X)||_F^2 &=& 
||\V_+(\X)^{1/2} \bold{\Gamma}_+(\X) \V_+(\X)^{1/2}||_F^2 
+ ||\V_-(\X)^{1/2} \bold{\Gamma}_-(\X) \V_-(\X)^{1/2}||_F^2 \\
& & + 2 \gamma_{\max}(\X)^2 ||\P_+(\X) \X \P_-(\X)||_F^2 \\
& \le & ||\V_+^{1/4}(\X) \bold{\Gamma}_+(\X) \V_+(\X)^{1/4}||_F^2 
+ ||\V_-^{1/4}(\X) \bold{\Gamma}_-(\X) \V_-(\X)^{1/4}||_F^2 
+ 2 \gamma_{\max}(\X)^2 \frac{n^3}{4} \\
& \le & N(\X) + \frac{1}{2} M_1^2 n^3.
\end{eqnarray*}
\hfill$\square$

We put Lemma~\ref{le:Dub} into Lemma~\ref{le:progress} to obtain a new upper bound on
$q(\alpha_k, \X^k)$;
\begin{eqnarray}
q(\alpha_k, \X^k) 
&\le& f(\X^k) 
- \frac{1}{2} \min \left\{\frac{N(\X^k)^2}{M_2 \left(N(\X^k) + \frac{1}{2} M_1^2 n^3 \right) },
\frac{N(\X^k)}{M_1},
\frac{\Delta_k N(\X^k)}{\sqrt{N(\X^k) + \frac{1}{2} M_1^2 n^3}}  \right\}. \label{eq:fupdate2}
\end{eqnarray}

In Algorithm~\ref{al:trust}, we call the $k$th iteration a {\it successful} iteration if 
$\X^{k+1}$ is set as $\overline{\X}^k$ in Step~4,
that is, $r_k \ge \mu_1$. Otherwise, the $k$th iteration is called an {\it unsuccessful} iteration.
For a successful iteration, we obtain a decrease in the objective function 
\begin{eqnarray}
f(\X^{k+1}) & \le &  f(\X^k) - \mu_1 (f(\X^k) - q(\alpha_k, \X^k)) 
\nonumber \\
&\le & f(\X^k)
- \frac{\mu_1}{2} \min \left\{\frac{N(\X^k)^2}{M_2 \left( N(\X^k) + \frac{1}{2} M_1^2 n^3 \right)},
\frac{N(\X^k)}{M_1},
\frac{\Delta_k N(\X^k)}{\sqrt{N(\X^k) + \frac{1}{2} M_1^2 n^3}}  \right\}. \label{eq:fupdate}
\end{eqnarray}
Since it holds $f(\X^{k+1}) = f(\X^k)$ for an unsuccessful iteration,
the objective value $f(\X^k)$ is non-increasing in Algorithm~\ref{al:trust}.

We are now prepared to show that
there exists a subsequence of $\{N(\X^k)\}$ that converges to zero.

\begin{THEO} \label{th:sub}
When $\{\X^k\}$ is the sequence generated by Algorithm~\ref{al:trust} with 
the stopping threshold $\epsilon = 0$,
it holds that 
\begin{eqnarray*}
\liminf_{k \to \infty} N(\X^k) = 0.
\end{eqnarray*}
\end{THEO}

\noindent{\bf Proof:}
We assume that there exists $\hat{\epsilon} > 0$ and $k_0$ such that 
$N(\X^k) \ge \hat{\epsilon}$ for any $k \ge k_0$, and we will derive a contradiction.

Let $\KC = \{k_1, k_2, \ldots, k_i, \ldots \} $ be the successful iterations. 
If $\KC$ is a finite sequence, let $k_i$ be the last iteration of $\KC$.
Since all of the iterations after $k_i$ are unsuccessful, the update rule of $\Delta_{k}$ 
(Step 5 of Algorithm~\ref{al:trust}) implies
$\Delta_{k_i + j} = \eta_1^j \Delta_{k_i}$.
Hence, we obtain $\lim_{j \to \infty} \Delta_j = 0$.
Next, we consider the case when $\KC$ is an infinite sequence.
The function $\frac{x^2}{x + \frac{1}{2} M_1^2 n^3}$ is an increasing function 
for $x > 0$, 
so that it holds from (\ref{eq:fupdate}) that for $k_i \in \KC$,
\begin{eqnarray*}
 f(\X^{k_i+1}) 
&\le& f(\X^{k_i}) 
- \frac{\mu_1}{2} \min 
\left\{\frac{N(\X^{k_i})^2}{M_2 \left(N(\X^{k_i}) + \frac{1}{2} M_1^2 n^3 \right) },
\frac{N(\X^{k_i})}{M_1},
\frac{\Delta_{k_i} N(\X^{k_i})}{\sqrt{N(\X^{k_i}) + \frac{1}{2} M_1^2 n^3}}  \right\} \\
&\le &
f(\X^{k_i}) 
- \frac{\mu_1}{2} \min 
\left\{\frac{\hat{\epsilon}^2}{M_2 \left(\hat{\epsilon} + \frac{1}{2} M_1^2 n^3 \right) },
\frac{\hat{\epsilon}}{M_1},
\frac{\Delta_{k_i} \hat{\epsilon}}{\sqrt{\hat{\epsilon} + \frac{1}{2} M_1^2 n^3}}  \right\}.
\end{eqnarray*}
Since $f$ is continuous on a closed set $\FC$ and $\X^k \in \FC$ for each $k$,
$f(\X^{k_i})$ is bounded below.
Therefore, it holds $\lim_{i \to \infty} \Delta_{k_i} = 0$.
From 
Step 5 of Algorithm~\ref{al:trust}, 
 it holds that 
$\Delta_{j}  \le \eta_2 \Delta_{k_{i}}$
for the unsuccessful iterations $j = k_{i}+1 \ldots, k_{i+1} - 1$.
Hence, we obtain $\lim_{j \to \infty} \Delta_j = 0$, regardless of the finiteness of $\KC$.
From (\ref{eq:fupdate2}) and
$N(\X^k) \ge \hat{\epsilon}$, 
it holds for sufficiently large $k$ that 
\begin{eqnarray}
f(\X^k)  - q(\alpha_k,\X^k) 
\ge \frac{1}{2} \frac{\Delta_k \hat{\epsilon}}{\sqrt{\hat{\epsilon} + \frac{1}{2} M_1^2 n^3}}
> 0. \label{eq:denominator}
\end{eqnarray}

We will take a close look at the ratio $r_k$.
From the Taylor expansion, there exists $\xi \in [0,1]$ such that 
{\small
\begin{eqnarray*}
f(\X^k - \alpha_k \S(\X^k)) &=& f(\X^k) - \alpha_k \langle \nabla f(\X^k) \ | \ \S(\X^k) \rangle
+ \frac{\alpha_k^2}{2} \langle \S(\X^k) \
| \ \nabla^2 f(\X^k -  \xi \alpha_k \S(\X^k)) \ | \ \S(\X^k) \rangle.
\end{eqnarray*}
} 
Therefore, 
{\small
\begin{eqnarray*}
|f(\overline{\X}^k) - q(\alpha^k, \X^k)| 
&\le& \frac{\alpha_k^2}{2} 
\left| \langle \S(\X^k) \
| \ \nabla^2 f(\X^k -  \xi \alpha_k \S(\X^k)) \ | \ \S(\X^k) \rangle 
- \langle \S(\X^k) \
| \ \nabla^2 f(\X^k) \ | \ \S(\X^k) \rangle \right|  \\
&\le& \frac{\Delta_k^2}{2} (M_2 + M_2) = \Delta_k^2 M_2.
\end{eqnarray*}
} 

Using (\ref{eq:denominator}) in the denominator, the ratio $r_k$ is evaluated by
\begin{eqnarray*}
|r_k - 1| = \frac{|f(\overline{\X}^k) - q(\alpha_k, \X^k)| }{|f(\X^k) - q(\alpha_k, \X^k)|}  
\le \frac{\Delta_k^2 M_2}{\frac{1}{2} \frac{\Delta_k \hat{\epsilon}}{\sqrt{\hat{\epsilon} + \frac{1}{2} M_1^2 n^3}}} 
= \Delta_k \frac{2 M_2 \sqrt{\hat{\epsilon} + \frac{1}{2} M_1^2 n^3}}{\hat{\epsilon}}.
\end{eqnarray*}
Therefore, $\lim_{k \to \infty} \Delta_k = 0$ leads to  $\lim_{k \to \infty} r_k = 1$.
From Step 5 of Algorithm~\ref{al:trust}, 
we have $\Delta_{k + 1} \ge \Delta_{k} $ for sufficiently large $k$.
Thus, there exists $\hat{k}_0$ such that $\Delta_{k} \ge \Delta_{\hat{k}_0}$ 
for $\forall k \ge \hat{k}_0$,
but  this contradicts
$\lim_{k \to \infty} \Delta_k = 0$.
Hence, $\liminf_{k \to \infty} N(\X^k) = 0$.
\hfill$\square$

\subsection{Convergence of the whole sequence} \label{sec:whole}

Using the convergence of the subsequence, we will show in Theorem~\ref{th:whole} 
that the whole sequence of $\{N(\X^k)\}$
converges to zero. 
We use the following two lemmas to prove Theorem~\ref{th:whole}.

\begin{LEMM}\label{le:Mub}
For $\X \in \FC$ and $\A, \B \in \SMAT^n$, we have
\begin{eqnarray*}
\left| \langle \nabla f(\X) \ | \ \A \rangle \right | 
&\le& \sqrt{n} M_1 ||\A||_F \\
\left| \langle \A \ | \ \nabla^2 f(\X) \ | \ \B \rangle \right | 
&\le& 3 M_2 ||\A||_F ||\B||_F.
\end{eqnarray*}
\end{LEMM}

\noindent{\bf Proof:}
The first inequality holds by
$\left| \langle \nabla f(\X) \ | \ \A \rangle \right | \le ||\nabla f(\X)||_F ||\A||_F$  
and 
$||\nabla f(\X)||_F \le \sqrt{n}||\nabla f(\X)||_2$
from \cite[(1.2.27)]{YUAN06}.

For the second inequality, we start with the following inequality derived from
the definition of $M_2$;
\begin{eqnarray*}
\left| \langle \D \ | \ \nabla^2 f(\X) \ | \ \D \rangle \right|
\le M_2 ||\D||_F^2 \quad \mbox{for} \quad \forall \D \in \SMAT^n.
\end{eqnarray*}
Therefore, we get 
$\left |\langle \A \ | \ \ \nabla^2 f(\X) \ | \ \A \rangle \right| \le M_2 ||\A||_F^2$
and
$\left |\langle \B \ | \ \nabla^2 f(\X) \ | \ \B \rangle\right| \le M_2 ||\B||_F^2$.
Furthermore, we put $\A - t \B$ into $\D$
to obtain the following inequality, which holds for any $t \in \Real$;
\begin{eqnarray*}
\left| \langle \A - t \B \ | \ \nabla^2 f(\X) \ | \ \A - t \B \rangle \right|
\le M_2 ||\A - t \B ||_F^2. 
\end{eqnarray*}
Therefore, the inequality 
\begin{eqnarray*}
(M_2 ||\B||_F^2 - \langle \B \ | \ \nabla^2 f(\X) \ | \ \B \rangle ) t^2 
- 2 (M_2 \langle \A \ | \ \B \rangle - \langle \A \ | \ \nabla^2 f(\X) \ | \ \B \rangle ) t & & \\
+ (M_2 ||\A||_F^2 - \langle \A \ | \ \nabla^2 f(\X) \ | \ \A \rangle ) &\ge& 0
\end{eqnarray*}
holds for any $t \in \Real$, and we can derive
\begin{eqnarray*}
& & (M_2 \langle \A \ | \ \B \rangle - \langle \A \ | \ \nabla^2 f(\X) | \ \B \rangle )^2 \\
&\le& \left(M_2 ||\A||_F^2 - \langle \A \ | \ \nabla^2 f(\X) \ | \ \A \rangle \right)
\left( M_2 ||\B||_F^2 - \langle \B \ | \ \nabla^2 f(\X) \ | \ \B \rangle \right) \\
&\le& (2 M_2 ||\A||_F^2) (2 M_2 ||\B||_F^2). 
\end{eqnarray*}
Consequently, it holds that 
\begin{eqnarray*}
\langle \A \ | \ \nabla^2 f(\X) \ | \ \B \rangle 
&\le& M_2 \langle \A \ | \ \B \rangle + \sqrt{(2 M_2 ||\A||_F^2) (2 M_2 ||\B||_F^2)}  \\
&\le& M_2 ||\A||_F ||\B||_F + 2 M_2 ||\A||_F ||\B||_F
= 3 M_2 ||\A||_F ||\B||_F 
\end{eqnarray*}
In addition, we replace $\A$ with $-\A$ to obtain
\begin{eqnarray*}
\langle -\A \ | \ \nabla^2 f(\X) \ | \ \B \rangle 
\le 3 M_2 ||\A||_F ||\B||_F.
\end{eqnarray*}
By combining these inequalities, we get 
$\left| \langle \A \ | \ \nabla^2 f(\X) \ | \ \B \rangle \right | 
\le 3 M_2 ||\A||_F ||\B||_F$.

\hfill$\square$

\begin{LEMM}\label{le:fbelow}
For $\X^k \in \FC$, it holds that $\underline{f}(\X^k) \ge - n \sqrt{ N(\X^k)}.$
\end{LEMM}
\noindent{\bf Proof:}

The objective function of 
(\ref{eq:underlinef}) at $\X \in \FC$ can be evaluated from below by
{\small
\begin{eqnarray*}
 \langle \nabla f(\X^k) \ | \ \X - \X^k \rangle 
& = & \langle \P_+(\X^k) \bold{\Gamma}_+(\X^k) \P_+(\X^k)^T 
+ \P_-(\X^k) \bold{\Gamma}_-(\X^k) \P_-(\X^k)^T 
 \ | \ \X - \X^k \rangle \\
& = & 
  \langle \bold{\Gamma}_+(\X^k) \ | \ \P_+(\X^k)^T \X \P_+(\X^k) \rangle 
- \langle  \bold{\Gamma}_-(\X^k) \ | \ \P_-(\X^k)^T (\I - \X) \P_-(\X^k) \rangle \\
& & - \langle \bold{\Gamma}_+(\X^k) \ | \ \V_+(\X^k) \rangle 
+ \langle \bold{\Gamma}_-(\X^k) \ | \ \V_-(\X^k) \rangle \\
& \ge & 
- \langle \bold{\Gamma}_+(\X^k) \ | \ \V_+(\X^k) \rangle 
+ \langle \bold{\Gamma}_-(\X^k) \ | \ \V_-(\X^k) \rangle.
\end{eqnarray*}
}
Furthermore, an upper bound of
$\langle \bold{\Gamma}_+(\X^k) \ | \ \V_+(\X^k) \rangle$ is given by
\begin{eqnarray*}
\langle \bold{\Gamma}_+(\X^k) \ | \ \V_+(\X^k) \rangle
&=& 
 Trace(\V_+(\X^k)^{1/4} \V_+(\X^k)^{1/4} \bold{\Gamma}_+(\X^k)
   \V_+(\X^k)^{1/4} \V_+(\X^k)^{1/4}) \\
&  \le & ||\V_+(\X^k)^{1/4}||_F 
|| \V_+(\X^k)^{1/4} \bold{\Gamma}_+(\X^k) \V_+(\X^k)^{1/4} ||_F 
||\V_+(\X^k)^{1/4}||_F  \\
&\le& n_+(\X^k) || \V_+(\X^k)^{1/4} \bold{\Gamma}_+(\X^k) \V_+(\X^k)^{1/4} ||_F.
\end{eqnarray*}
Here, we used $||\V_+(\X^k)^{1/4}||_F \le n_+(\X^k)$ derived from $\O \preceq \V_+(\X^k)^{1/4} \preceq \I$.
In a similar way, it also holds 
$\langle - \bold{\Gamma}_-(\X^k) \ | \ \V_-(\X^k) \rangle
\le  n_-(\X^k) || \V_-(\X^k)^{1/4} \bold{\Gamma}_-(\X^k) \V_-(\X^k)^{1/4} ||_F$.

Therefore, we obtain
\begin{eqnarray*}
\underline{f}(\X^k) &\ge &
- n_+(\X^k) || \V_+(\X^k)^{1/4} \bold{\Gamma}_+(\X^k) \V_+(\X^k)^{1/4} ||_F
- n_-(\X^k) || \V_-(\X^k)^{1/4} \bold{\Gamma}_-(\X^k) \V_-(\X^k)^{1/4} ||_F \\
&\ge& - (n_+(\X^k) + n_-(\X^k)) \\
& & \times
\sqrt{|| \V_+(\X^k)^{1/4} \bold{\Gamma}_+(\X^k) \V_+(\X^k)^{1/4} ||_F^2
 + || \V_-(\X^k)^{1/4} \bold{\Gamma}_-(\X^k) \V_-(\X^k)^{1/4} ||_F^2} \\
&=& -n \sqrt{N(\X^k)}.
\end{eqnarray*}

For the second inequality, we used an inequality
$a b + c d \le (a+c) \sqrt{b^2 + d^2}$ for $a,b,c,d \ge 0$.


\hfill$\square$

We are ready to prove the convergence of the whole sequence.

\begin{THEO} \label{th:whole}
When $\{\X^k\}$ is the sequence generated by Algorithm~\ref{al:trust} with $\epsilon = 0$,
it holds that 
\begin{eqnarray*}
\lim_{k \to \infty} N(\X^k) = 0.
\end{eqnarray*}
\end{THEO}
\noindent{\bf Proof:}

We take a small positive number $\epsilon_1$ such that $ 0 < \epsilon_1 \le 16 n^2 M_1^2$.
We assume that there is an infinite subsequence
$\KC := \{k_1, k_2, \ldots, k_i, \ldots  \} \subset \{1, 2, \ldots \}$ such that 
$N(\X^{k_i}) \ge \epsilon_1$ for $\forall k_i \in \KC$, and we will derive a contradiction.

From Theorem~\ref{th:sub}, we can take a subsequence 
$\LC := \{l_1, l_2, \ldots, l_i, \ldots  \} \subset \{1, 2, \ldots \}$  such that 
\begin{eqnarray*}
\left\{\begin{array}{lccl}
N(\X^{k}) & \ge &  \epsilon_2^2  & \quad \mbox{for} \quad k = k_i, k_i+1, \ldots, l_i -1\\
N(\X^{l_i}) & < & \epsilon_2^2.  & \\
\end{array}\right.
\end{eqnarray*}
where $\epsilon_2 := \frac{\epsilon_1}{4 n M_1}$.
Note that this is consistent with $N(\X^{k_i}) \ge \epsilon_1$, 
since we took $ 0 < \epsilon_1 \le 16 n^2 M_1^2$.

If the $k$th iteration is a successful iteration and $k_i \le k < l_i$, 
we put $N(\X^k) \ge \epsilon_2^2$ into 
(\ref{eq:fupdate}) and obtain
\begin{eqnarray*}
f(\X^{k+1}) 
&\le & f(\X^k) 
- \frac{\mu_1}{2} \min \left\{\frac{\epsilon_2^4}{M_2 \left( \epsilon_2^2 + \frac{1}{2} M_1^2 n^3 \right)},
\frac{\epsilon_2^2}{M_1},
\frac{\Delta_k \epsilon_2^2}{\sqrt{\epsilon_2^2 + \frac{1}{2} M_1^2 n^3}}  \right\}. 
\end{eqnarray*}
Since $f$ is bounded below, 
if $k$ is sufficiently large, it holds that 
\begin{eqnarray*}
f(\X^{k+1}) 
&\le & f(\X^k) 
 - \Delta_k \epsilon_3
\end{eqnarray*}
where $\epsilon_3 := 
\frac{\mu_1}{2} 
\frac{\epsilon_2^2}{\sqrt{\epsilon_2^2 + \frac{1}{2} M_1^2 n^3}}$.
We update the matrix with $\X^{k+1} = \X^{k} - \alpha_k \S(\X^k)$ in a successful iteration, 
therefore, 
we use $\alpha_k \le \Delta_k$ and $||\S(\X^k)||_F = \frac{||\D(\X^k)||_F}{||\D(\X^k)||_F} = 1$ to derive
\begin{eqnarray*}
||\X^{k} - \X^{k+1}||_F  \le \Delta_k
\le \frac{f(\X^k) - f(\X^{k+1})}{\epsilon_3}.
\end{eqnarray*}
This inequality is also valid when the $k$th iteration is an unsuccessful iteration, 
since the matrix is updated with $\X^{k+1} = \X^{k}$.
Hence, it holds that 
\begin{eqnarray*}
& & ||\X^{k_i} - \X^{l_i}||_F \\
&\le & ||\X^{k_i} - \X^{k_i+1}||_F 
+ ||\X^{k_i+1} - \X^{k_i+2}||_F \ldots
+ ||\X^{l_i-1} - \X^{l_i}||_F \\
&\le& \frac{1}{\epsilon_3}
\left( (f(\X^{k_i}) - f(\X^{k_i+1})) + 
(f(\X^{k_i+1}) - f(\X^{k_i+2})) + \cdots + (f(\X^{l_i-1}) - f(\X^{l_i})) \right) \\
&=& \frac{f(\X^{k_i}) -f(\X^{l_i})}{\epsilon_3}.
\end{eqnarray*}
Since the objective function $f(\X^k)$ is non-increasing and bounded below,
this implies that $\lim_{i \to \infty} ||\X^{k_i} - \X^{l_i}||_F = 0$. Therefore,
for $\epsilon_4 := \frac{\sqrt{n} \epsilon_2}{M_1 + 3 M_2} > 0$,
there exists $i_0$ such that 
$||\X^{k_i} - \X^{l_i}||_F < \epsilon_4$ for $\forall i \ge i_0$.


Since $\X^k \in \FC$, it holds that $- \I \preceq \X - \X^{k_i} \preceq \I$ for $\X \in \FC$.
Therefore, we have 
an inequality
$||\X - \X^{k_i}||_F \le \sqrt{n}$.
For $\X \in \FC$ and $i \ge i_0$, 
it holds that 
\begin{eqnarray*}
& & \left| \langle \nabla f(\X^{k_i})  \ | \ \X - \X^{k_i} \rangle  
- \langle \nabla f(\X^{l_i})  \ | \ \X - \X^{l_i} \rangle  \right| \\
& = & \left| \langle \nabla f(\X^{l_i} + (\X^{k_i} - \X^{l_i})) \ | \ \X - \X^{k_i} \rangle  
- \langle \nabla f(\X^{l_i})  \ | \ \X - \X^{l_i} \rangle  \right| \\
& = & \left| 
\langle \nabla f(\X^{l_i}) \ | \ \X - \X^{k_i} \rangle  
+ \int_0^1 \langle \X^{k_i} - \X^{l_i} \ | \ 
\nabla^2 f(\X^{l_i} + \xi(\X^{k_i} - \X^{l_i})) 
\ | \ \X - \X^{k_i} \rangle d\xi \right. \\
& & \hspace{1cm}
\left. - \langle \nabla f(\X^{l_i})  \ | \ \X - \X^{l_i} \rangle  \right| \\
& = & 
\left| \int_0^1 
\langle \X^{k_i} - \X^{l_i} \ | \ \nabla^2 f(\X^{l_i} + \xi(\X^{k_i} - \X^{l_i})) 
\ | \ \X - \X^{k_i} \rangle  d\xi
 - \langle \nabla f(\X^{l_i})  \ | \ \X^{k_i} - \X^{l_i} \rangle  \right| \\
& \le & 
\int_0^1 \left|
\langle \X^{k_i} - \X^{l_i} \ | \ \nabla^2 f(\X^{l_i} + \xi(\X^{k_i} - \X^{l_i})) 
\ | \ \X - \X^{k_i} \rangle  \right|  d\xi
 + \left| \langle \nabla f(\X^{l_i})  \ | \ \X^{k_i} - \X^{l_i} \rangle  \right| \\
&\le& 3 M_2 ||\X^{k_i} - \X^{l_i}||_F
|| \X - \X^{k_i}||_F
+ \sqrt{n} M_1 ||\X^{k_i} - \X^{l_i}||_F \\
&\le&  3 M_2  \epsilon_4  \sqrt{n} + \sqrt{n} M_1 \epsilon_4
 = \sqrt{n}(3 M_2 + M_1) \epsilon_4
= n \epsilon_2.
\end{eqnarray*}
Here, we used Lemma~\ref{le:Mub} for the second inequality.
Hence, we have
\begin{eqnarray}
 \langle \nabla f(\X^{k_i})  \ | \ \X - \X^{k_i} \rangle  
\ge \langle \nabla f(\X^{l_i})  \ | \ \X - \X^{l_i} \rangle  
- n \epsilon_2.\label{eq:fkl}
\end{eqnarray}

If $\gamma_{\max}(\X^{k_i}) = 0$, then $\nabla f(\X^{k_i}) = \O$,
and this results in $N(\X^{k_i}) = 0$ from (\ref{eq:N}).
Therefore, from the assumption $N(\X^{k_i}) \ge \epsilon_2^2$
we know that $\gamma_{\max}(\X^{k_i}) > 0$.
Since $\X^{k_i} - \frac{\D(\X^{k_i})}{\gamma_{\max}(\X^{k_i})} \in \FC$ from 
Lemma~\ref{le:rangealpha}, 
we can put $\X^{k_i} - \frac{\D(\X^{k_i})}{\gamma_{\max}(\X^{k_i})}$
into (\ref{eq:fkl}) to get
\begin{eqnarray*}
 \langle \nabla f (\X^{k_i}) \ | \  -\frac{\D(\X^{k_i})}{\gamma_{\max}(\X^{k_i})} \rangle 
\ge \langle \nabla f (\X^{l_i}) \ 
| \ \left( \X^{k_i} -\frac{\D(\X^{k_i})}{\gamma_{\max}(\X^{k_i})} \right) - \X^{l_j} \rangle
- n \epsilon_2
\ge \underline{f}(\X^{l_j}) - n \epsilon_2.
\end{eqnarray*}

With Lemma~\ref{le:fbelow} and $N(\X^{l_i}) < \epsilon_2^2$, 
we have an upper bound on $N(\X^{k_i})$;
\begin{eqnarray*}
N(\X^{k_i}) &=&
\langle \nabla f (\X^{k_i}) \ | \  \D(\X^{k_i})  \rangle  
\le \gamma_{\max}(\X^{k_i}) 
(-\underline{f}(\X^{l_j})  + n \epsilon_2)  \\
&\le& \gamma_{\max}(\X^{k_i}) 
(n \sqrt{N(\X^{l_j})}  + n \epsilon_2) 
\le M_1 (n \epsilon_2 + n \epsilon_2) = 2 n M_1  \epsilon_2. 
\end{eqnarray*}

Therefore, we obtain the contradiction;
\begin{eqnarray*}
\epsilon_1 \le N(\X^{k_i}) \le 2 n M_1 \epsilon_2
= \frac{1}{2} \epsilon_1
< \epsilon_1.
\end{eqnarray*}
Hence, $\lim_{k \to \infty} N(\X^k) = 0.$

\hfill$\square$

Combining Lemma~\ref{le:fbelow} and Theorem~\ref{th:whole}, 
we derive the property for the first-order optimality condition.

\begin{THEO}\label{th:flimit}
When $\{\X^k\}$ is the sequence generated by Algorithm~\ref{al:trust} with $\epsilon = 0$,
it holds that 
\begin{eqnarray*}
\lim_{k \to \infty} \underline{f}(\X^{k}) = 0.
\end{eqnarray*}
\end{THEO}
\noindent{\bf Proof:}
From Lemma~\ref{le:fbelow}, 
we know that 
$-n \sqrt{N(\X^k)} \le \underline{f}(\X^k) \le 0$. 
Hence, Theorem~\ref{th:whole} indicates
$\lim_{k \to \infty} \underline{f}(\X^k) = 0$.

\hfill$\square$

Using Theorem~\ref{th:flimit}, we can show an additional result on the convergence.
To make the generated sequence $\{\X^k\}$ itself converge, we need a stronger assumption
on the objective function,
for example, strong convexity.

\begin{CORO}\label{co:strongly-convex}
If the objective function $f$ is strongly convex, that is, there exists $\nu > 0$ such that 
\begin{eqnarray*}
f(\Y) \ge f(\X) + \langle \nabla f(\X) \ | \ \Y - \X \rangle
+ \frac{\nu}{2} ||\Y - \X||_F^2 \quad \mbox{for} \quad
\forall \X, \forall \Y \in \FC,
\end{eqnarray*}
then the sequence $\{\X^k\}$ generated by Algorithm~\ref{al:trust} with $\epsilon = 0$ converges.
Furthermore, the accumulation point $\X^* := \lim_{k \to \infty} \X^k$ 
is an optimal solution.
\end{CORO}
\noindent{\bf Proof:}

From $\X^k \in \FC$ and the definition of $\underline{f}(\X^j)$ for $\X^j \in \FC$,
we have an inequality
$\underline{f}(\X^j) \le \langle \nabla f(\X^j) \ | \ \X^k - \X^j \rangle$.
By swapping $\X^k$ and $\X^j$, we also obtain the inequality
$\underline{f}(\X^k) \le \langle \nabla f(\X^k) \ | \ \X^j - \X^k \rangle$.
The addition of these two inequalities results in
\begin{eqnarray*}
\langle \nabla f(\X^k) - \nabla f(\X^j) \ | \ \X^k - \X^j \rangle
\le - \underline{f}(\X^k) - \underline{f}(\X^j).
\end{eqnarray*}
Theorem 2.1.9 of \cite{NESTEROV04} gives equivalent conditions of strong convexity, and 
one of them is 
\begin{eqnarray*}
\langle \nabla f(\Y) - \nabla f(\X) \ | \ \Y - \X \rangle 
\ge \nu ||\Y - \X||_F^2 \quad \forall \X, \forall \Y \in \FC.
\end{eqnarray*}
Due to this inequality, we get
\begin{eqnarray*}
||\X^k - \X^j||_F \le \frac{1}{\nu} 
\sqrt{- \underline{f}(\X^k) - \underline{f}(\X^j)}.
\end{eqnarray*}
Theorem~\ref{th:flimit} implies that the sequence $\{\X^k\}$ is a Cauchy
sequence. Since $\{\X^k\}$ is generated in the closed and bounded set $\FC$,
it converges to a point of $\FC$. 
Hence, the accumulation point $\X^* = \lim_{k \to \infty} \X^k$
satisfies the first-order optimality condition.
From the assumption that the objective function is convex,
$\X^*$ is an optimal solution.
\hfill$\square$

\section{Numerical Results}\label{sec:results}

To evaluate the performance of the proposed method,
we conducted a numerical test.
The computing environment was Debian Linux run on  AMD Opteron Processor 4386 (3 GHz) 
and 128 GB of memory space, and we used Matlab R2014a.

The test functions used are listed below and they are classified into the two groups.
The functions of Group I were selected from \cite{XU11}, and
we added new functions as Group II. Function 5 and 6 are 
an extension of Generalized Rosenbrock function~\cite{NASH84} and its variant
with cosine functions, respectively.

\begin{enumerate}[Group I:]
\item 
\begin{enumerate}[{Function} 1.]
\item $f(\X) = -2 \langle \C_1 \ | \  \X \rangle + \langle \X \ | \  \X \rangle$;
\item $f(\X) = 3 \cos (\langle \X \ | \ \X \rangle) + \sin (\langle \X + \C_1 \ | \ \X + \C_1 \rangle)$;
\item $f(\X) = \log(\langle \X \ | \ \X \rangle + 1) + 5 \langle \C_1 \ | \ \X \rangle$;
\end{enumerate}
\item 
\begin{enumerate}[{Function} 1.]
\setcounter{enumii}{3}
\item $f(\X) = 1 + 2 \frac{\langle \X - \C_1 \ |  \X - \C_1 \rangle^3}{n^3}$;
\item $f(\X) =  1
     + \sum_{i=1}^n \sum_{j=i}^n (A_{ij} - X_{ij})^2$ \\
     \hspace{1cm} $+ 100 \sum_{i=1}^{n-1}\sum_{j=i}^{n-1}
\left(\frac{A_{ij}^2}{A_{i,j+1}} X_{i,j+1}-X_{ij}^2\right)^2 $\\
     \hspace{1cm} $+ 100 \sum_{i=1}^{n-1}\left(\frac{A_{in}^2}{A_{i+1,i+1}}X_{i+1,i+1}-X_{i,n}^2 \right)^2$;
\item $ f(\X) = 
\frac{1}{n^2}  \sum_{i=1}^{n} 
\left( \sum_{j=1, j \ne i}^n \frac{X_{ij}}{A_{ij}} - (n-1) \frac{X_{ii}^2}{A_{ii}^2}  \right)^2$ 
\\
   \hspace{1cm} 
$ - \frac{1}{n^2} \sum_{i=1}^{n} \sum_{j=1}^{n} \cos ((X_{ij}-A_{ij})^2)$;
\item $f(\X) = \langle \C_1 \ | \  \X \rangle - \log \det (\X + \bar{\epsilon} \I) 
- \log \det ((1+\bar{\epsilon})\I -\X)$;
\end{enumerate}
\end{enumerate}

To generate the matrix $\C_1$ in Functions 1, 4, and 7, we chose the eigenvalues 
$\kappa_1,\ldots, \kappa_n$ randomly from the interval $[-1,2]$ and 
multiply a randomly-generated orthogonal matrix $\Q$, namely, 
$\C_1 := \Q diag(\kappa_1, \ldots, \kappa_n) \Q^T$.
The elements $A_{ij}$ in Functions 5 and 6 were set as
$A_{ii} =  \frac{1}{2}$ for $i=1,\ldots,n$ and $A_{ij} = \frac{1}{2(n-1)}$ for $i \ne j$.
The parameter $\bar{\epsilon}$ in Function 7 was set as $\bar{\epsilon} = 0.02$.

We compared the performance of three methods, 
PIM  (the proposed iterative method, Algorithm~\ref{al:trust}), 
FEAS (the feasible direction method of Xu {\it et. al.}~\cite{XU11}),
and PEN (the penalty barrier method \cite{BENTAL97, KOCVARA03} implemented in PENLAB~\cite{PENLAB13}).
We started PIM and FEAS with the initial point $\X^0 = \frac{1}{2} \I$,
while PEN automatically chose its initial point.
and $\Delta^0 = 1$.
We used the following condition as the stopping criterion;
\begin{center}
\begin{tabular}{ll}
PIM  & $N(\X^k) < 10^{-7}$ or $\frac{|f(\X^k)-f(\X^{k-1})|}{\max\{|f(\X^k)|,1\}} < 10^{-6}$ \\
FEAS  & $|Trace(\bold{\Gamma}_{-}(\X^k)) - \langle f(\X^k) \ | \ \X^k \rangle| < 10^{-6}$ or $\frac{|f(\X^k)-f(\X^{k-1})|}{\max\{|f(\X^k)|,1\}} < 10^{-6}$ \\
PEN  & the default parameter of PENLAB.
\end{tabular}
\end{center}
For details of the stopping criterion on FEAS and PEN, refer to \cite{XU11} and \cite{PENLAB13}, respectively.
We also stopped the computation when the computation time exceeded 24 hours.

Tables~\ref{table:results1} and \ref{table:results2} show the numerical results
of Group I and Group II, respectively.  The first column is the function type,
and the second column $n$ is the size of the matrix $\X$.
The third column indicates the method we applied,
and the fourth column is the objective value.
The fifth column is the number of main iterations,
and the six column is the computation time in seconds.
The last three columns correspond to the evaluation count of
the function value $f(\X)$, the gradient matrix $\nabla f(\X)$, and 
the Hessian mapping $\nabla^2 f(\X)$.


\begin{table}[htbp]
\begin{center}
\caption{Numerical results on Group I.}
\label{table:results1}
{ \footnotesize
\begin{tabular}{|r|r|c|r|r|r|r|r|r|}
\hline
type & \multicolumn{1}{|c|}{$n$} & \multicolumn{1}{|c|}{method} 
& \multicolumn{1}{|c|}{obj} & \multicolumn{1}{|c|}{iter} 
& \multicolumn{1}{|c|}{cpu} & co.$f$ & co.$\nabla f$ & co.$\nabla^2 f$ \\
\hline
 1 & 50 &   PIM & $-3.631 \times 10^1$ & 48 & 0.08 & 95 & 48 & 48 \\
 1 & 50 & FEAS & $-3.633 \times 10^1$ & 36 & 0.04 & 215 & 36 & 0 \\
 1 & 50 &  PEN & $-3.633 \times 10^1$ & 22 & 323.70 & 62 & 31 & 22 \\
\hline
 1 & 100 &   PIM & $-7.930 \times 10^1$ & 67 & 0.30 & 133 & 67 & 67 \\
 1 & 100 & FEAS & $-7.932 \times 10^1$ & 36 & 0.11 & 239 & 36 & 0 \\
 1 & 100 &  PEN & $-7.932 \times 10^1$ & 23 & 5554.30 & 64 & 32 & 23 \\
\hline
 1 & 500 &   PIM & $-3.572 \times 10^2$ & 81 & 7.70 & 161 & 81 & 81 \\
 1 & 500 & FEAS & $-3.574 \times 10^2$ & 37 & 2.24 & 250 & 37 & 0 \\
\hline
 1 & 1000 &   PIM & $-8.648 \times 10^2$ & 64 & 30.16 & 127 & 64 & 64 \\
 1 & 1000 & FEAS & $-8.651 \times 10^2$ & 32 & 9.62 & 204 & 32 & 0 \\
\hline
 1 & 5000 &   PIM & $-3.861 \times 10^3$ & 80 & 3497.86 & 159 & 80 & 80 \\
 1 & 5000 & FEAS & $-3.862 \times 10^3$ & 36 & 1111.89 & 232 & 36 & 0 \\
\hline
 1 & 10000 &   PIM & $-7.731 \times 10^3$ & 73 & 24730.04 & 145 & 73 & 73 \\
 1 & 10000 & FEAS & $-7.734 \times 10^3$ & 34 & 7782.18 & 213 & 34 & 0 \\
\hline
\hline
 2 & 50 &   PIM & $-4.000$ & 23 & 0.04 & 45 & 23 & 23 \\
 2 & 50 & FEAS & $-4.000$ & 31 & 0.04 & 293 & 31 & 0 \\
 2 & 50 &  PEN & $-4.000$ & 115 & 1808.54 & 1857 & 124 & 116 \\
\hline
 2 & 100 &   PIM & $-4.000$ & 40 & 0.19 & 79 & 40 & 40 \\
 2 & 100 & FEAS & $-4.000$ & 13 & 0.05 & 122 & 13 & 0 \\
 2 & 100 &  PEN & $-3.985$ & 15 & 4581.35 & 114 & 21 & 18 \\
\hline
 2 & 500 &   PIM & $-4.000$ & 26 & 2.40 & 51 & 26 & 26 \\
 2 & 500 & FEAS & $-4.000$ & 17 & 1.31 & 183 & 17 & 0 \\
\hline
 2 & 1000 &   PIM & $-4.000$ & 17 & 6.32 & 33 & 17 & 17 \\
 2 & 1000 & FEAS & $-4.000$ & 10 & 4.30 & 134 & 10 & 0 \\
\hline
 2 & 5000 &   PIM & $-4.000$ & 28 & 1205.24 & 55 & 28 & 28 \\
 2 & 5000 & FEAS & $-4.000$ & 13 & 449.59 & 161 & 13 & 0 \\
\hline
 2 & 10000 &   PIM & $-4.000$ & 27 & 8461.01 & 53 & 27 & 27 \\
 2 & 10000 & FEAS & $-3.951$ & 8 & 2009.73 & 73 & 8 & 0 \\
\hline
\hline
 3 & 50 &   PIM & $-3.756 \times 10^1$  & 201 & 0.35 & 401 & 201 & 201 \\
 3 & 50 & FEAS & $-3.756 \times 10^1$ & 2 & 0.01 & 3 & 2 & 0 \\
 3 & 50 &  PEN & $-3.756 \times 10^1$  & 28 & 418.41 & 76 & 36 & 28 \\
\hline
 3 & 100 &   PIM & $-7.418 \times 10^1$ & 208 & 0.93 & 415 & 208 & 208 \\
 3 & 100 & FEAS & $-7.419 \times 10^1$ & 7 & 0.02 & 24 & 7 & 0 \\
 3 & 100 &  PEN &  $-7.419 \times 10^1$ & 30 & 7316.99 & 81 & 37 & 30 \\
\hline
 3 & 500 &   PIM & $-3.625 \times 10^2$ & 257 & 25.62 & 513 & 257 & 257 \\
 3 & 500 & FEAS & $-3.625 \times 10^2$ & 2 & 0.13 & 3 & 2 & 0 \\
\hline
 3 & 1000 &   PIM & $-7.739 \times 10^2$ & 269 & 128.23 & 537 & 269 & 269 \\
 3 & 1000 & FEAS & $-7.741 \times 10^2$ & 2 & 0.65 & 3 & 2 & 0 \\
\hline
 3 & 5000 &   PIM & $-4.129 \times 10^3$ & 257 & $11996.45$ & 513 & 257 & 257 \\
 3 & 5000 & FEAS & $-4.129 \times 10^3$ & 2 & 74.63 & 3 & 2 & 0 \\
\hline
 3 & 10000 &   PIM & $-8.294 \times 10^3$ & 256 & 92901.29 & 511 & 256 & 256 \\
 3 & 10000 & FEAS & $-8.295 \times 10^3$ & 2 & 575.84 & 3 & 2 & 0 \\
\hline
\end{tabular}
} 
\end{center}
\end{table}

\begin{table}[htbp]
\begin{center}
\caption{Numerical results on Group II.}
\label{table:results2}
{\footnotesize
\begin{tabular}{|r|r|c|r|r|r|r|r|r|}
\hline
type & \multicolumn{1}{|c|}{$n$} & \multicolumn{1}{|c|}{method} 
& \multicolumn{1}{|c|}{obj} & \multicolumn{1}{|c|}{iter} 
& \multicolumn{1}{|c|}{cpu} & co.$f$ & co.$\nabla f$ & co.$\nabla^2 f$ \\
\hline
 4 & 50 &   PIM & $1.041$ & 26 & 0.07 & 51 & 26 & 26 \\
 4 & 50 & FEAS & $1.041$ & 16 & 0.03 & 77 & 16 & 0 \\
 4  & 50 &  PEN & $1.041$ & 23 & 328.07 & 75 & 37 & 23 \\
\hline
 4 & 100 &   PIM & $1.039$ & 42 & 0.19 & 83 & 42 & 42 \\
 4 & 100 & FEAS & $1.039$ & 24 & 0.07 & 138 & 24 & 0 \\
 4 & 100 &  PEN & $1.039$ & 25 & 5824.27 & 83 & 40 & 25 \\
\hline
 4 & 500 &   PIM & $1.024$ & 21 & 2.08 & 41 & 21 & 21 \\
 4 & 500 &  FEAS & $1.024$ & 23 & 1.35 & 123 & 23 & 0 \\
\hline
 4 & 1000 &   PIM & $1.023$ & 14 & 6.77 & 27 & 14 & 14 \\
 4 & 1000 & FEAS & $1.023$ & 25 & 7.35 & 142 & 25 & 0 \\
\hline
 4 & 5000 &   PIM & $1.024$ & 12 & 517.62 & 23 & 12 & 12 \\
 4  & 5000 & FEAS & $1.024$ & 25 & 715.68 & 134 & 25 & 0 \\
\hline
 4 & 10000 &   PIM & $1.025$ & 12 & 4140.26 & 23 & 12 & 12 \\
 4 & 10000 & FEAS & $1.025$ & 21 & 4866.45 & 109 & 21 & 0 \\
\hline
\hline
 5 & 50 &   PIM & $1.122$ & 4 & 0.01 & 7 & 4 & 4 \\
 5 & 50 & FEAS & $1.126$ & 19 & 0.05 & 252 & 19 & 0 \\
 5 & 50 &  PEN & $1.000$ & 20 & 294.58 & 61 & 30 & 20 \\
\hline
 5 & 100 &   PIM & $1.117$ & 6 & 0.06 & 11 & 6 & 6 \\
 5 & 100 & FEAS & $1.125$ & 16 & 0.15 & 226 & 16 & 0 \\
 5 & 100 &  PEN & $1.000$  & 20 & 4814.74 & 61 & 30 & 20 \\
\hline
 5 & 500 &   PIM & $1.004$ & 4 & 0.82 & 7 & 4 & 4 \\
 5 & 500 & FEAS & $1.125$ & 16 & 4.28 & 286 & 16 & 0 \\
\hline
 5 & 1000 &   PIM & $1.008$ & 4 & 4.02 & 7 & 4 & 4 \\
 5 & 1000 & FEAS & $1.125$ & 18 & 26.96 & 352 & 18 & 0 \\
\hline
 5 & 5000 &   PIM & $1.002$ & 4 & 192.25 & 7 & 4 & 4 \\
 5 & 5000 & FEAS & $1.125$ & 90 & 6345.17 & 2279 & 90 & 0 \\
\hline
 5 & 10000 &   PIM & $1.013$ & 4 & 1332.14 & 7 & 4 & 4 \\
 5 & 10000 & FEAS & $1.124$ & 122 & 51611.04 & 3285 & 122 & 0 \\
\hline
\hline
6 & 50 &   PIM & $-1.000$ & 20 & 0.11 & 39 & 20 & 20 \\
6 & 50 & FEAS & $-1.000$ & 12 & 0.10 & 92 & 12 & 0 \\
6 & 50 &  PEN & $-1.000$ & 300 & 4577.01 & 915 & 1218 & 300 \\
\hline
6 & 100 &   PIM & $-1.000$ & 20 & 0.36 & 39 & 20 & 20 \\
6 & 100 & FEAS & $-1.000$ & 16 & 0.56 & 150 & 16 & 0 \\
6 & 100 &  PEN & $-9.997 \times 10^{-1}$ & 300 & 73262.02 & 1005 & 1308 & 300 \\
\hline
6 & 500 &   PIM & $-1.000$ & 18 & 10.00 & 35 & 18 & 18 \\
6 & 500 & FEAS & $-1.000$ & 12 & 9.36 & 110 & 12 & 0 \\
\hline
6 & 1000 &   PIM & $-1.000$ & 4 & 9.42 & 7 & 4 & 4 \\
6 & 1000 & FEAS & $-1.000$ & 12 & 56.33 & 110 & 12 & 0 \\
\hline
6 & 5000 &   PIM & $-1.000$ & 4 & 406.01 & 7 & 4 & 4 \\
6 & 5000 & FEAS & $-1.000$ & 13 & 2046.55 & 130 & 13 & 0 \\
\hline
6 & 10000 &   PIM & $-1.000$ & 3 & 2076.17 & 5 & 3 & 3 \\
6 & 10000 & FEAS & $-1.000$ & 14 & 10416.77 & 130 & 14 & 0 \\
\hline
\hline
7 & 50 &   PIM & $7.817 \times 10^1$ & 10 & 0.03 & 19 & 10 & 10 \\
7 & 50 & FEAS & $7.817 \times 10^1$ & 15 & 0.06 & 108 & 15 & 0 \\
7 & 50 &  PEN & $7.817 \times 10^1$ & 13 & 195.41 & 38 & 19 & 13 \\
\hline
7 & 100 &   PIM & $1.583 \times 10^2$ & 10 & 0.11 & 19 & 10 & 10 \\
7 & 100 & FEAS & $1.583 \times 10^2$ & 17 & 0.30 & 13 & 17 & 0 \\
7 & 100 &  PEN & $1.583 \times 10^2$ & 14 & 3427.86 & 40 & 20 & 14 \\
\hline
7 & 500 &   PIM & $7.825 \times 10^2$ & 10 & 2.73 & 19 & 10 & 10 \\
7 & 500 & FEAS & $7.825 \times 10^2$ & 16 & 6.22 & 116 & 16 & 0 \\
\hline
7 & 1000 &   PIM & $1.556 \times 10^3$ & 10 & 12.02 & 19 & 10 & 10 \\
7 & 1000 & FEAS & $1.556 \times 10^3$ & 10 & 15.79 & 60 & 10 & 0 \\
\hline
7 & 5000 &   PIM & $7.707 \times 10^3$ & 11 & 1708.11 & 21 & 11 & 11 \\
7 & 5000 & FEAS & $7.707 \times 10^3$ & 16 & 4931.70 & 115 & 16 & 0 \\
\hline
7 & 10000 &   PIM & $1.533 \times 10^4$ & 11 & 13379.96 & 21 & 11 & 11 \\
7 & 10000 & FEAS & $1.533 \times 10^4$ & 14 & 32643.49 & 94 & 14 & 0 \\
\hline
\end{tabular}
} 
\end{center}
\end{table}

From these tables, PEN was much slow compared to PIM and FEAS.
We did not include the results of PEN for large problems $ n \ge 500$,
since PEN did not finish the computation for $ n = 500$ in 24 hours.
Though it attained better solution for Function 5,
PENLAB~\cite{PENLAB13} handled the symmetric matrix $\X$ 
as $n(n+1)/2$ independent variables ($X_{11}$, $X_{12}$, $\ldots$, $X_{1n}$,
$X_{22}$, $\ldots$, $X_{2n}$, $\ldots$, $X_{nn}$),
and it stored all the elements of the Hessian mapping $\nabla^2 f(\X)$, therefore,
the computation cost  was estimated as $\OC(n^4)$ from \cite{KOCVARA03}.
This heavy cost restricted PENLAB to the small sizes.
PIM also used the information of the Hessian mapping,
but in only the scalar value $\langle \S \ | \ \nabla^2 f(\X) \ | \ \S \rangle$.
Hence, the computation cost of each iteration in PIM is much lower than PEN,
and this low cost is the key to handling large problems.


In the comparison between PIM and FEAS, the computation time of FEAS
was shorter than PIM in Table~\ref{table:results1}, but longer in Table~\ref{table:results2}.
The functions in Group I involved the variable matrix $\X$ in the linear form 
$\langle \C_1 \ | \ \X \rangle $ or the quadratic form $\langle \X \ | \ \X \rangle$,
and this simple structure was favorable for the feasible direction method, which was based on 
a steepest descent direction.
In contrast, the functions in Group II have stronger nonlinearity than Group I.
The evaluation count with respect to the function value (co.$f$) implies that 
this stronger nonlinearity demanded FEAS have a large number of back-step loop.
In particular, FEAS needed many iterations for Rosenbrock-type functions (Functions~5 and 6). 
PIM reduced the number of iterations by 
the properties of the search direction $D(\X)$
and the quadratic approximation with the Hessian mapping.
In particular, $\D(\X)$ encompassed the information of the distance to
the boundary to the box-constraints as $\bold{\V}_+(\X)$ and $\bold{\V}_-(\X)$.
Therefore, PIM was faster than FEAS for the functions of Group II.

\section{Conclusions and Future Directions}\label{sec:conclusions}

In this paper, we proposed an iterative method for box-constrained SDPs.
The search direction $\D(\X)$ studied in Section~2 enabled us 
to include the information of the distance from the current point to the boundary of the feasible set $\FC$.
We discussed the convergence property 
of the generated sequence. 
The numerical tests in Section~4 showed that 
the proposed method was more favorable 
for functions with strong nonlinearity
than the feasible direction method,
mainly due to the distance information included in $\D(\X)$.
In addition, the proposed method handled the larger problems than the penalty barrier method,
since our method did not hold the Hessian mapping in memory space.

One of future researches
would be the combination of the feasible direction and the proposed method,
since the feasible direction method fits simple functions.
For such a combination, 
we should extend the convergence analysis from this paper.
Another point is the convergence for a second-order optimality condition,
as proven in \cite{COLEMAN96} for box-constrained problem~(\ref{eq:boxLP}).
The proof in~\cite{COLEMAN96} required further stronger assumptions than this paper and
the second-order optimality condition for nonlinear semidefinite programs
involves not only the Hessian mapping but also an additional mapping~\cite{SHAPRIO97},
so we remain it as a matter to be discussed further.


\bibliography{reference}

\end{document}